\documentclass[review,3p,11pt]{elsarticle}
\biboptions{sort&compress}

\usepackage{amsmath,amssymb,amsthm}
\usepackage{bm,color}
\usepackage{stmaryrd}
\usepackage[T1]{fontenc}
\usepackage{float, caption, subcaption}
\usepackage{booktabs}
\usepackage{multirow}
\usepackage[bookmarks=true,
  bookmarksnumbered=true,
  bookmarksopen=true,
  colorlinks,
  pdfborder=001,
  linkcolor=black]{hyperref}
\usepackage[ruled,linesnumbered]{algorithm2e}

\allowdisplaybreaks

\newtheorem{remark}{Remark}[section]

\numberwithin{equation}{section}
\numberwithin{figure}{section}
\numberwithin{table}{section}
\newcommand\bbR{\mathbb{R}}

\newcommand\bq{\bm{q}}
\newcommand\bv{\bm{v}}
\newcommand\bV{\bm{V}}
\newcommand\bu{\bm{u}}

\newcommand\bx{\bm{x}}

\newcommand\dd{\mathrm{d}}

\newcommand\norm[1]{\left\lVert #1 \right\rVert}
\newcommand\diag{\mathrm{diag}}

\newcommand\encop{\mathfrak{E}}
\newcommand\decop{\mathfrak{D}}
\newcommand\enc{{\tt enc}}
\newcommand\dec{{\tt dec}}

\newcommand\rb{\tt{rb}}

\journal{Journal of Computational Physics}

\begin{document}

\begin{frontmatter}



\title{Non-intrusive data-driven reduced-order modeling for time-dependent parametrized problems}


\author[inst1]{Junming Duan\corref{cor1}}
\affiliation[inst1]{organization={Chair of Computational Mathematics and Simulation Science, \'Ecole Polytechnique F\'ed\'erale de Lausanne},
            city={Lausanne},
            postcode={1015},
            country={Switzerland}}
\cortext[cor1]{Corresponding author}
\ead{junming.duan@epfl.ch}

\author[inst1]{Jan S. Hesthaven}
\ead{jan.hesthaven@epfl.ch}

\begin{abstract}
Reduced-order models are indispensable for multi-query or real-time problems.
However, there are still many challenges to constructing efficient ROMs for time-dependent parametrized problems.
Using a linear reduced space is inefficient for time-dependent nonlinear problems, especially for transport-dominated problems.
The non-linearity usually needs to be addressed by hyper-reduction techniques,
such as DEIM, but it is intrusive and relies on the assumption of affine dependence of parameters.
This paper proposes and studies a non-intrusive reduced-order modeling approach for time-dependent parametrized problems.
It is purely data-driven and naturally split into offline and online stages.
During the offline stage,
a convolutional autoencoder, consisting of an encoder and a decoder, is trained to perform dimensionality reduction.
The encoder compresses the full-order solution snapshots to a nonlinear manifold or a low-dimensional reduced/latent space.
The decoder allows the recovery of the full-order solution from the latent space.
To deal with the time-dependent problems, a high-order dynamic mode decomposition (HODMD)
is utilized to model the trajectories in the latent space for each parameter.
During the online stage,
the HODMD models are first utilized to obtain the latent variables at a new time,
then interpolation techniques are adopted to recover the latent variables at a new parameter value,
and the full-order solution is recovered by the decoder.
Some numerical tests are conducted to show that the approach can be used to
predict the unseen full-order solution at new times and parameter values fast and accurately,
including transport-dominated problems.
\end{abstract}




\end{frontmatter}


\section{Introduction}\label{section:Intro}
In the study of many real-world applications in science and engineering,
such as optimization, control, or uncertainty quantification,
many time-dependent problems are described as parametrized partial differential equations (PDEs).
The parameters may come from physical properties, geometric configurations, initial or boundary conditions, etc.
Such parametrized PDEs can be solved by standard numerical methods,
e.g., the finite difference method, finite volume method, spectral method, finite element method, etc.
However, thousands of degrees of freedom are usually required to obtain sufficiently accurate solutions,
i.e., high-fidelity solutions of the full-order model (FOM),
which leads to high demands on computational resources.
Furthermore, in the context of multi-query or real-time tasks,
these PDEs need to be solved for a large number of different parameter values.
In such cases, the development of efficient low-dimensional models
that allow a fast evaluation of an output of interest at a new time and parameter value with controlled loss of accuracy is of great interest.

During the past decades, the reduced-basis model (RBM) \cite{Hesthaven2015Certified, Quarteroni2016Reduced} has been developed to tackle this issue.
The key idea of the RBM is to replace the FOM with a surrogate model
by finding low-dimensional structures in a collection of full-order solutions at sampling times and parameter values,
called the snapshots, which describe the underlying spatial-temporal dynamics of the solution manifold.
One of the most popular methods is the projection-based RBM \cite{Benner2015A},
where a linear combination of bases spans a low-dimensional approximation (reduced subspace) of the solution manifold,
and then the FOM is projected into the low-dimensional reduced subspace to obtain the RBM.
To recover the optimal linear low-dimensional approximation,
the proper orthogonal decomposition (POD) is utilized
to give the dominant orthonormal modes by decomposing the snapshot matrix based on the singular value decomposition (SVD).
The projection-based method splits into offline and online stages,
and most of the computational costs are completed during the offline stage, including the collection of the snapshots
and construction of the RBM,
which can be deployed to recover a fast response to the request at the online stage,
since the computational cost scales with the dimension of the RBM.

The development of the RBM has been well studied for the elliptic, stationary, and linear problems with certified error control.
The readers are referred to the books \cite{Hesthaven2015Certified, Quarteroni2016Reduced},
the recent review article \cite{Hesthaven2022Reduced}, and references therein.
Additional challenges come from the non-linearity in the FOM,
resulting in the computational cost of the nonlinear terms in the RBM scaling with the high dimension of the FOM.
In such cases, hyper-reduction strategies \cite{Ryckelynck2009Hyper} are used.
Most of these methods rely on sparse sampling through interpolation of the nonlinear operators,
e.g., missing point estimation \cite{Astrid2008Missing},
the empirical interpolation method (EIM) \cite{Barrault2004An},
the discrete empirical interpolation method (DEIM) \cite{Chaturantabut2010Nonlinear},
Gauss–Newton with approximated tensors (GNAT) \cite{Carlberg2013The}, etc.
However, these methods are generally intrusive, need access to the original full-order solvers,
and an efficient implementation may be non-trivial.
This motivates the design of non-intrusive methods.
In \cite{Hesthaven2018Non}, a neural network was used to approximate the map
from the parameter space to the reduced coefficients, and shown to be efficient.
Gaussian process regression has been also used to build non-intrusive RBM in \cite{Guo2018Reduced}.

Another difficulty is that the dimension of the linear reduced space generated by the POD can be very high,
in cases when the Kolmogorov $n$-width decays very slowly,
e.g., a moving-front solution to the advection equation \cite{Ehrlacher2020nonlinear}.
In this case, the efficiency is limited by the high dimensionality of the RBM.
In \cite{Carlberg2015Adaptive}, $h$-adaptivity was used to enrich the reduced bases.
The quadratic operator inference \cite{Peherstorfer2016Data} was proposed to build a non-intrusive projection-based RBM.
Neural networks are also utilized to construct low-dimensional reduced space,
e.g., an autoencoder used in \cite{Lee2020Model, Maulik2021Reduced},
which was shown to be superior to the classic POD-based linear subspace.

In this paper, we are concerned with time-dependent problems,
thus the RBM should capture the underlying dynamics.
The adaptive reduced bases and sampling via low-rank updates were proposed in \cite{Peherstorfer2015Online, Peherstorfer2020Model},
and can be viewed as $r$-adaptivity.
The recurrent neural network (RNN), more specifically, long short-term memory (LSTM) network,
has also been adopted to model the dynamics in \cite{Maulik2021Reduced},
but only the states at some discrete times are available, not the whole trajectory.
Another way is to view the time as another parameter and build the map from the time and parameter to the reduced coefficients
using a neural network \cite{Fresca2021Deep, Fresca2022POD}.

In this work, we will build the surrogate model for the latent dynamics based on
the higher-order extension of the dynamic mode decomposition (DMD), i.e. HODMD.
Such methods have been used for the model reduction of time-dependent problems,
and are interpretable through Koopman spectral theory \cite{Brunton2022Modern}.
The key idea is to choose suitable measurements/coordinates,
and transform the original nonlinear dynamical system to a linear one acting on these new coordinates.
In the DMD, the full-order solutions are directly chosen as the coordinates,
while in the HODMD, the time-delay embedding is introduced to enrich the coordinates.
The first step of the original DMD algorithms is to perform dimensionality reduction based on the POD,
while we will use the latent variables from the nonlinear dimensionality reduction as input,
leading to more efficient low-dimensional representation.

This paper proposes and studies a new non-intrusive data-driven reduced-order modeling approach for time-dependent parametrized problems,
split into the offline and online stages.
During the offline stage, the reduced/latent space is generated by a deep convolutional autoencoder,
which has been shown to result in an approximation that is more efficient than the linear subspace.
The trained encoder yields the latent variables for the snapshots,
and the HODMD is employed to build the surrogate models for the latent dynamics at each training parameter value.
During the online stage, for a given new time and parameter value,
the HODMD models are first utilized to obtain the latent variables at the new time,
then the latent variables at the new parameter value are obtained by interpolation,
and the full-order solution is recovered by the decoder.
Our approach is purely data-driven and does not use a priori knowledge of the underlying physical model.
Three tests are conducted to verify our method,
i.e., 1D Burgers' equation, 2D Rayleigh-B\'enard convection, and 2D Kelvin-Helmholtz instability.
The results show that the approach can work well for transport-dominated problems with a low-dimensional latent space.
It can be used to predict the unseen full-order solution at new times and parameter values fast and accurately.
It is loosely coupled so that one can adjust each component and tune the corresponding parameters,
making it easy to control the errors in each part.

This paper is organized as follows.
Section \ref{section:FOM} introduces the FOM and the collection of the snapshots.
Our data-driven non-intrusive method will be detailed in Section \ref{section:NumerMethod}.
Some numerical results are presented in Section \ref{section:Result} to validate the effectiveness and performance of our method,
with concluding remarks in Section \ref{section:Conc}.

\section{Full-order model formulation}\label{section:FOM}
This paper considers a set of general parametrized partial differential equations (PDEs)
\begin{equation}\label{eq:FOM}
	\left\{
	\begin{aligned}
		&\dfrac{\partial }{\partial t} \bu(t, \bx; \bm{\omega}) + \bm{f}(\bu, \nabla\bu, \nabla^2\bu, \cdots; \bm{\omega}) = 0, \\
		&\bu(t=0, \bx; \bm{\omega}) = \bu_0(\bx; \bm{\omega}),
	\end{aligned}
	\right.
\end{equation}
where $\bu(t, \bx; \bm{\omega})\in\bbR^{m}$ is the solution vector with $m$ components at time $t\in[0,+\infty)$,
spatial coordinate $\bx\in\bbR^{d}$, depending on the parameter vector $\bm{\omega}\in\bbR^{d_p}$,
and $\bm{f}$ is a nonlinear function involving the spatial derivatives of $\bu$.
The model \eqref{eq:FOM} can be solved by using the finite difference, finite volume, finite element, or spectral methods
for a given parameter value $\bm{\omega}$ to get the following semi-discrete full-order model (FOM)
\begin{equation}\label{eq:DFOM}
	\left\{
	\begin{aligned}
		&\dfrac{\dd }{\dd t} \bu_h(t; \bm{\omega}) + \bm{f}_h(\bu_h; \bm{\omega}) = 0, \\
		&\bu_h(t=0; \bm{\omega}) = \bu_{h,0}(\bm{\omega}),
	\end{aligned}
	\right.
\end{equation}
where $\bm{f}_h$ is the discrete spatial operator.
By integrating \eqref{eq:DFOM} in time one obtains the full-order solution
$\bu_h(t; \bm{\omega})\in\bbR^{n_h}$ at time $t$.
In this paper, $n_h=m\times n_y\times n_x$,
where $n_x, n_y$ are the numbers of spatial degrees of freedom in the $x$- and $y$-direction, respectively
($n_y=1$ for 1D problem),
although our approach can be extended to three dimensions without any difficulty.

\section{Non-intrusive data-driven method}\label{section:NumerMethod}
This section will introduce some existing results,
and then present our reduced-order modeling approach for the discrete FOM \eqref{eq:DFOM}.

\subsection{Linear model order reduction: POD-Galerkin approach}\label{sec:PODG}
The classic projection-based RBMs are widely used in many applications \cite{Benner2015A, Hesthaven2015Certified, Quarteroni2016Reduced}.
In such methods, the solution manifold is approximated by a reduced $n_{\rb}$-dimensional linear space $\mathcal{S}_{n_{\rb}}$
spanned by the $n_{\rb}$ columns of a matrix $\mathcal{V}_{\rb}\in\bbR^{n_h\times n_{\rb}}$
\begin{equation}\label{eq:rb_app}
	\bu_h(t; \bm{\omega}) = \mathcal{V}_{\rb}\bu_{\rb}(t; \bm{\omega}) + \overline{\bu_h},
\end{equation}
where $\bu_{\rb}$ are the coefficients of the reduced basis function,
and $\overline{\bu_h}$ is a reference solution.
The POD is one of the most popular methods to recover the orthonormal modes, or $\mathcal{V}_{\rb}$,
that span the best linear subspace based on the singular value decomposition (SVD).
If one collects snapshots at the sampling times $t_1,t_2,\cdots,t_{n_t}$ and parameter values $\bm{\omega}_1,\bm{\omega}_2,\cdots,\bm{\omega}_{n_p}$
and form the snapshot matrix
\begin{equation*}
	S=[\bu_h(t_1;\bm{\omega}_1)-\overline{\bu_h}, \bu_h(t_{n_t};\bm{\omega}_1)-\overline{\bu_h}, \dots, \bu_h(t_{n_t};\bm{\omega}_{n_p})-\overline{\bu_h}]\in\bbR^{n_h\times n_s}
\end{equation*}
the SVD decomposition of the snapshot matrix is
\begin{equation*}
	S = \mathcal{\bV} \Sigma \widetilde{\mathcal{\bV}}^\mathrm{T},
\end{equation*}
where $\Sigma = \text{diag}\left\{\sigma_{1}, \sigma_{2}, \cdots, \sigma_{r} \right\}\in \bbR^{n_h\times n_s}$
with the singular values $\sigma_{1}\geqslant \sigma_{2}\geqslant \cdots\geqslant \sigma_{r}\geqslant 0$ and $r= \min\{n_h,n_s\}$.
Here the columns of $\mathcal{\bV}\in\bbR^{n_h\times n_h}$ and $\widetilde{\mathcal{\bV}}\in\bbR^{n_s\times n_s}$ are the orthonormal left and right singular vectors.
Then the matrix $\mathcal{V}_{\rb}$ is selected as the first $n_{\rb}$ columns of $\mathcal{V}$,
and the dimension $n_{\rb}$ can be determined by the relative energy threshold
\begin{equation}\label{eq:SVD_epsilon}
	{\sum\limits_{l=1}^{n_{\rb}} \sigma_l^2}\Big/{\sum\limits_{l=1}^{r} \sigma_l^2} \geqslant 1-\epsilon,
\end{equation}
with $\epsilon$ a small number close to zero.
By the Schmidt-Eckart-Young theorem \cite{Eckart1936The, Schmidt1907On},
the matrix $\mathcal{V}_{\rb}$ minimizes the projection error
\begin{equation*}
	\sum_{i=1}^{n_t}\sum_{j=1}^{n_p}\norm{\bu_h(t_i;\bm{\omega}_j) -
	\mathcal{V}_{\rb}\mathcal{V}_{\rb}^\mathrm{T}\bu_h(t_i;\bm{\omega}_j)}^2_{\bbR^{n_h}}
	=\min_{\mathcal{W}\in\mathcal{S}_{n_{\rb}}} \sum_{i=1}^{n_t}\sum_{j=1}^{n_p}
	\norm{\bu_h(t_i;\bm{\omega}_j) - \mathcal{W}\mathcal{W}^\mathrm{T}\bu_h(t_i;\bm{\omega}_j)}^2_{\bbR^{n_h}}
\end{equation*}
over the set $\mathcal{S}_{n_{\rb}} = \{\mathcal{W}\in\bbR^{n_h\times n_{\rb}}: \mathcal{W}^\mathrm{T}\mathcal{W} = I_{n_{\rb}} \}$,
i.e., all the rank $n_{\rb}$ orthonormal bases.

Inserting the ansatz \eqref{eq:rb_app} into the semi-discrete FOM \eqref{eq:DFOM}, one recovers the projected RBM
\begin{equation}\label{eq:PODG}
	\left\{
	\begin{aligned}
		&\dfrac{\dd }{\dd t} \bu_{\rb}(t; \bm{\omega})
		+ \mathcal{\bV}_{\rb}^\mathrm{T}\bm{f}_h(\mathcal{\bV}_{\rb}\bu_{\rb} + \overline{\bu_h}) = 0, \\
		&\bu_{\rb}(t=0; \bm{\omega}) = \mathcal{\bV}_{\rb}^\mathrm{T}\left(\bu_{h,0}(\bm{\omega})-\overline{\bu_h}\right).
	\end{aligned}
	\right.
\end{equation}
Although the POD-Galerkin methods have been widely used,
the dimension of the linear reduced space can be large to obtain an accurate approximation of the snapshot matrix
for time-dependent nonlinear problems.
The Kolmogorov $n$-width \cite{Kolmogoroff1936Uber} provides one way to quantify
the approximation of the optimal $n_{\rb}$-dimensional linear trial subspace $\mathcal{S}_{n_{\rb}}$
to the discrete solution manifold $\mathcal{M}_h$, defined as
\begin{equation*}
	d_{n_{\rb}}(\mathcal{M}_h) = \inf_{\mathcal{S}_{n_{\rb}}}\sup_{f\in\mathcal{M}_h}\inf_{g\in\mathcal{M}_h}\norm{f-g}.
\end{equation*}
In some cases, e.g. for the elliptic problems of high regularity in the parameter space \cite{Hesthaven2015Certified},
the Kolmogorov $n$-width decays exponentially,
$d_{n_{\rb}}(\mathcal{M}_h)\leqslant Ce^{-cn_{\rb}}$.
However, for transport-dominated problems, e.g. moving front solutions to the linear advection equation,
the linear subspace is not efficient, $d_{n_{\rb}}(\mathcal{M}_h)\leqslant Cn_{\rb}^{-1/2}$ \cite{Ehrlacher2020nonlinear}.

Another difficulty arises when applying the RBM \eqref{eq:PODG} for nonlinear problems,
which requires the use of $\mathcal{\bV}_{\rb}\bu_{\rb}$, and the computational cost scales with the dimension of the FOM, which is large.
Hyper-reduction techniques must be adopted to deal with the nonlinear terms,
such as (discrete) empirical interpolation methods (EIM/DEIM) \cite{Barrault2004An, Chaturantabut2010Nonlinear},
which aim at recovering the affine dependence of the parameter.
However, for a general application, such an assumption of affine dependence of the parameter may not hold,
and the approximation cost may be high.

\subsection{Autoencoder based nonlinear latent manifold}\label{sec:CAE}
The difficulties mentioned in the last section motivate the use of nonlinear trial space for the reduced spaces,
and nonlinear approximation for the nonlinear terms.
The neural network has been shown as a powerful tool for function approximation.
In this work, we adopt an autoencoder to perform the dimensionality reduction.

\begin{figure}
		\centering
		\includegraphics[width=\linewidth]{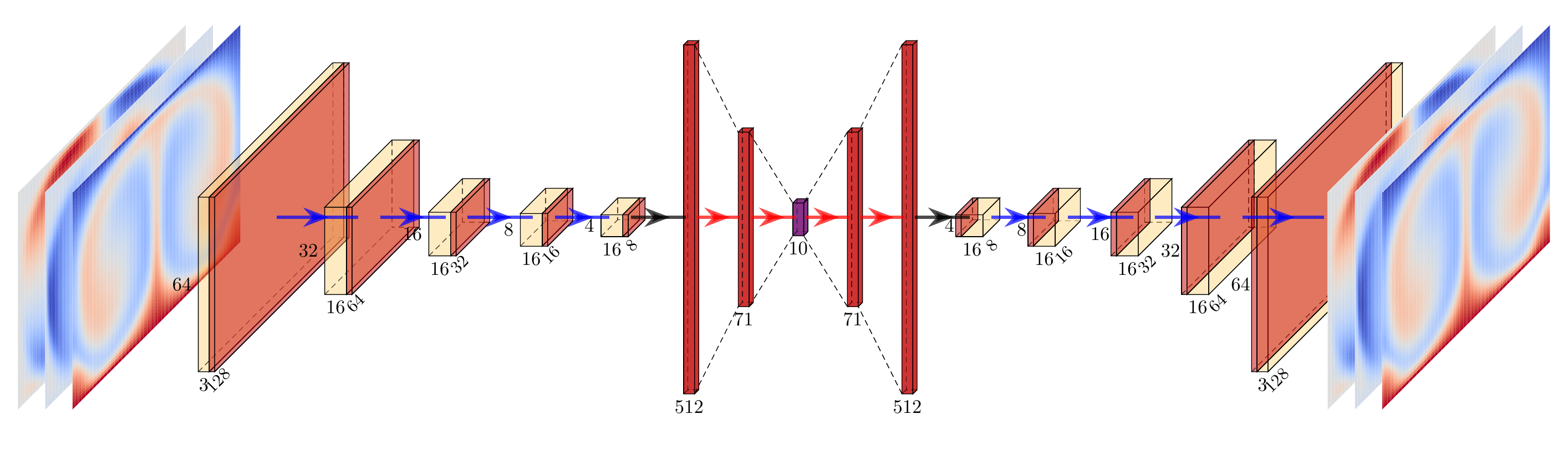}
		\caption{Architecture of the convolutional autoencoder used for 2D Rayleigh-B\'enard convection in Section \ref{sec:2DRBC}.}
		\label{fig:cae2d}
\end{figure}

The autoencoder consists of an encoder and a decoder part,
expressed as $\encop(\ \cdot\ ; \Theta_{\enc})$ and $\decop(\ \cdot\ ; \Theta_{\dec})$, respectively,
with $\Theta_{\enc}$ and $\Theta_{\dec}$ being the parameters (weights and biases) in the neural networks, recovered during the training.
Both are composed of multiple layers,
and the former compresses the full-order solutions to low-dimensional latent variables,
while the latter reconstructs the full-order solutions from the latent variables.
If only fully-connected layers are used, the autoencoder tends to have an extremely large number of parameters,
if the size of the input data is large, especially when considering solutions of the 2D or 3D FOMs.
In such a case, the training requires lots of data and can be very expensive.
As the full-order solutions on the structured meshes can be viewed as images,
we consider the convolutional neural networks (CNNs),
which have been shown to be very successful for compression in many image-related tasks.
CNN is very efficient since the total number of parameters is reduced due to parameter sharing.
This paper employs the convolutional autoencoder (CAE), also used in \cite{Lee2020Model, Maulik2021Reduced},
to produce a low-dimensional latent manifold.
A typical architecture of the convolutional autoencoder is shown in Figure \ref{fig:cae2d}.
Given an input $\bu_h\in\bbR^{m\times n_y\times n_x}$, the encoder part consists of multiple stacking layers as
\begin{equation*}
	\encop(\ \cdot\ ; \Theta_{\enc}) = \bm{h}_{n_L}(\ \cdot\ ; \Theta_{\enc, n_L}) \circ 
	\sigma_{n_L-1}\left(\bm{h}_{n_L-1}(\ \cdot\ ; \Theta_{\enc, n_L-1})\right) \circ \cdots \sigma_1\left(\bm{h}_{1}(\ \cdot\ ; \Theta_{\enc, 1})\right),
\end{equation*}
where the $i$th layer $\bm{h}_{i}, i=1,\cdots,n_L$ with the weights and biases $\Theta_{\enc, i}$
can be a convolutional layer or a fully-connected layer,
and $\sigma_i(\ \cdot\ )$ is the corresponding nonlinear activation function applied component-wise,
e.g., rectified linear unit (ReLU), hyperbolic tangent (Tanh), sigmoid linear unit (SiLU), etc.
We will use the SiLU function defined as
\begin{equation*}
	{\tt SiLU}(x) = x\ {\tt sigmoid}(x) = \dfrac{x}{1 + \exp(-x)}.
\end{equation*}
The convolutional layers reduce the spatial dimensions of the full-order solution
and change the number of channels.
The output after multiple convolutional layers is reshaped as a long vector,
and then transformed to the low-dimensional latent vector after several fully-connected layers.
For the decoder part, it is symmetric to the encoder part,
\begin{equation*}
	\decop(\ \cdot\ ; \Theta_{\dec}) = \bar{\bm{h}}_{1}(\ \cdot\ ; \Theta_{\dec, 1}) \circ 
	\sigma_{2}\left(\bar{\bm{h}}_{2}(\ \cdot\ ; \Theta_{\dec, 2})\right) \circ \cdots \sigma_{n_L}\left(\bar{\bm{h}}_{n_L}(\ \cdot\ ; \Theta_{\dec, n_L})\right),
\end{equation*}
where $\bar{\bm{h}}_{i}$ is a deconvolutional layer or a fully-connected layer,
and its input and output dimensions are opposite to those in $\bm{h}_{i}$ in the encoder.
The specific definition of the fully-connected layer, convolutional layer,
and deconvolutional layer can be found in the documents of PyTorch \cite{PyTorch2019}.
The $n_{\tt latent}$ dimensional central latent variables encode the main features and physical structures,
and they span a low-dimensional space, which is usually nonlinear.
If we only use one fully-connected layer without activation functions in both the encoder and the decoder,
the autoencoder is similar to the POD method,
and $n_{\tt latent}$ is just the number of the dominant modes chosen in the truncated SVD.

The loss function of the autoencoder is the reconstruction error defined as 
\begin{equation}\label{eq:cae_loss}
	\mathcal{L} = \dfrac{1}{n_tn_p}\sum\limits_{i=1}^{n_t}\sum\limits_{j=1}^{n_p}
	{\norm{\bu_h(t_i, \bm{\omega}_j)
				- \decop(\encop(\bu_h(t_i, \bm{\omega}_j); \Theta_{\enc}); \Theta_{\dec})}_2^2},
\end{equation}
where $\bu_h(t_i, \bm{\omega}_j)$ is the full-order solution at the time $t_i$ and parameter value $\bm{\omega}_j$,
and the summation is performed over the training set.
During the training, a gradient-based method is utilized to find the best parameters
$\Theta_{\enc}, \Theta_{\dec} = \mathop{\arg\min}\limits_{\Theta_{\enc}, \Theta_{\dec}} \mathcal{L}$.
After the training of the autoencoder, we use the encoder to obtain the latent variables $\bm{q} = \encop(\bu_h; \Theta_{\enc})$,
which span the latent manifold.

\begin{remark}\rm
	To avoid over-fitting, it is useful to add regularization terms in the loss function \eqref{eq:cae_loss},
	such as $\lambda\left(\norm{\Theta_{\enc}}_1 + \norm{\Theta_{\dec}}_1\right)$
	or $\lambda\left(\norm{\Theta_{\enc}}_2^2 + \norm{\Theta_{\dec}}_2^2\right)$,
	with $\lambda$ as the regularization parameter. 
	The $L^2$ regularization is equivalent to the weight decay technique in the standard stochastic gradient descent (SGD) optimizer,
	where the regularization terms appear directly in the update of the parameters.
	We use the weight decay technique implemented in the ADAM optimizer in PyTorch as it is found to be robust.
\end{remark}

\begin{remark}\rm
	In the numerical tests, we also employ early stopping technique, in other words,
	the training will be terminated if the error on the validation set has not become smaller than the best historical one
	for a given number of epochs. 
\end{remark}

\begin{remark}\rm
	In \cite{Fresca2021Deep, Fresca2022POD}, the convolutional autoencoder is used to generate the low-dimensional latent space,
	but the full-order solution is first reshaped as an image since it is not obtained on a structured mesh.
	The reshaping operator may destroy the spatial correlations, thus only full-order solutions on the structured meshes are considered in this paper.
\end{remark}

\subsection{HODMD for the temporal modeling of latent trajectory}\label{sec:HODMD}
To enable fast prediction at a new time, we build a surrogate model valid for the whole trajectory of the latent variables.
One way to model the temporal dynamics is to use neural networks,
such as recurrent neural networks (RNNs) \cite{Maulik2021Reduced}.
However, the RNNs can only give the states at discrete times, not the whole trajectory.
Furthermore, the error is known to accumulate fast in RNNs.
In this paper, we adopt the DMD \cite{Brunton2022Modern, Schmid2010Dynamic},
which is data-driven and has been proposed for the model reduction of time-dependent problems.
The classic DMD utilizes the POD to construct the reduced bases in space,
i.e., the DMD employs a linear subspace which may need a large number of  basis as mentioned in Section \ref{sec:PODG}.
We propose to reduce the degree of freedom in space by using the autoencoder first to obtain the latent variables
and then construct the latent dynamics based on the DMD.

Given the latent variables $\bm{q} = \encop(\bu_h; \Theta_{\enc})$ for the parameters $\bm{\omega}$ at $t=t_1,\cdots, t_{n_t}$,
we collect them as the snapshot matrix $\mathcal{Q} = \left[\bm{q}(t_1,\bm{\omega}), \cdots \bm{q}(t_{n_t},\bm{\omega})\right]$.
Then one can approximate the nonlinear dynamics of the latent variables by using a linear dynamical system
$$\mathcal{Q}_2 = \mathcal{A}\mathcal{Q}_1,\quad 
\mathcal{Q}_1 = \left[\bm{q}(t_1), \cdots \bm{q}(t_{n_t-1})\right], \quad \mathcal{Q}_2 = \left[\bm{q}(t_2), \cdots \bm{q}(t_{n_t})\right],$$
where $\bm{\omega}$ is omitted.
Algorithm \ref{alg:dmd} is used to compute the DMD eigenvalues, modes, and amplitudes.

\begin{algorithm}[H]
	\SetAlgoLined
	\KwIn{${\mathcal{Q}}_1, {\mathcal{Q}}_2, \epsilon$}
	\KwOut{$\Lambda, \Theta$}
	Compute the truncated SVD: ${\mathcal{Q}}_1 = \mathcal{X}\Sigma\widetilde{\mathcal{X}}^\mathrm{T}$,
	with $\epsilon$ the relative energy threshold \eqref{eq:SVD_epsilon} to select the number $L$ of the dominant singular vectors\;
	Obtain the reduced DMD operator $\widetilde{\mathcal{A}} = \mathcal{X}^\mathrm{T}{\mathcal{Q}}_2\widetilde{\mathcal{X}}\Sigma^{-1}$\;
	Compute the eigendecomposition $\widetilde{\mathcal{A}}\mathcal{Y} = \mathcal{Y} \Lambda$
	to get the reduced DMD modes $\mathcal{Y}$ and the DMD eigenvalues $\Lambda = \text{diag}\left\{\lambda_1,\cdots,\lambda_{L}\right\}$\;
	The DMD modes of the DMD operator ${\mathcal{A}}$ is defined as $\Xi = \left[{\bm{\xi}}_1,\cdots,{\bm{\xi}}_{L}\right] = {\mathcal{Q}}_2 \widetilde{\mathcal{X}} \Sigma^{-1} \mathcal{Y}$\;
	Predict the latent variable at $t=\tilde{t}$ by
	$\bm{q} (\tilde t) =  \sum\limits_{l=1}^{L} a_l \bm{\xi}_l {\lambda_l}^{(\tilde{t}-t_1)/\Delta t}$,
	with $\bm{a}=(a_1,a_2,\cdots,a_L)^\mathrm{T}$ the DMD amplitudes.
	\caption{DMD}
	\label{alg:dmd}
\end{algorithm}

The DMD amplitudes can be computed by $\bm{a} = \Xi^\dag\bm{q}(t_1)$.
Alternatively, the optimal amplitudes proposed in \cite{Jovanovic2014Sparsity},
which minimizes the errors between the reconstruction and all the snapshots,
obtained by solving the following optimization problem
\begin{equation*}
	\mathop{\arg\min}\limits_{\bm{a}} \norm{{\mathcal{Q}}_1 - \Xi~\diag\{a_1,a_2,\cdots,a_L\}~V_{\text{and}}}_F^2,
\end{equation*}
with the Vandermonde matrix
\begin{equation*}
	V_{\text{and}}=
	\begin{bmatrix}
		1 & \lambda_1 & \cdots & \lambda_1^{L-1} \\
		1 & \lambda_2 & \cdots & \lambda_2^{L-1} \\
		\vdots & \vdots & \ddots & \vdots \\
		1 & \lambda_L & \cdots & \lambda_L^{L-1} \\
	\end{bmatrix},
\end{equation*}
can be considered.
The solution is computed by solving the linear system
\begin{equation}\label{eq:amp_opt}
	\left((\Xi^*\Xi)\circ( \overline{V_{\text{and}}V_{\text{and}}^*})\right)\bm{a} =
	\overline{\diag(V_{\text{and}} {\mathcal{Q}}_1^* \Xi)}.
\end{equation}

For the latent dynamics considered in this paper,
the DMD is not accurate enough to capture the whole dynamics.
Thus in this paper, we propose to adopt the HODMD \cite{LeClainche2017Higher},
which can be viewed as a superimposed DMD containing more information in a sliding window.
The so-called time-delay embedding is used in HODMD as
\begin{equation*}
	\widehat{\bm{q}}_k = [\bm{q}(t_k)^\mathrm{T}, \bm{q}(t_{k+1})^\mathrm{T},\cdots, \bm{q}(t_{k+n_{\tt delay}-1})^\mathrm{T}]^\mathrm{T},
\end{equation*}
and the snapshot matrix is formed as the Hankel matrix
\begin{equation*}
	\widehat{\mathcal{Q}}_{1} =
	\left[\widehat{\bm{q}}_1, \widehat{\bm{q}}_2, \cdots, \widehat{\bm{q}}_{n_t-n_{\tt delay}}\right],\quad
	\widehat{\mathcal{Q}}_{2} =
	\left[\widehat{\bm{q}}_2, \widehat{\bm{q}}_3, \cdots, \widehat{\bm{q}}_{n_t-n_{\tt delay}+1}\right].
\end{equation*}
Then the following problem is solved $$\widehat{\mathcal{Q}}_2 = \mathcal{A}\widehat{\mathcal{Q}}_1$$
by using Algorithm \ref{alg:dmd}, except that for the prediction in the fourth step,
we only take the first $n_{\tt latent}$ components to recover the $n_{\tt latent}$ dimensional latent variables.
	
\begin{remark}\rm
	The parameter $n_{\tt delay}$ should not be close to one as in such cases there is not enough time-delay embedding, leading to inaccurate results.
	One also cannot choose too large $n_{\tt delay}$, since there are only a small number of modified snapshots in $\widehat{\mathcal{Q}}_1, \widehat{\mathcal{Q}}_2$,
	and the results at the trailing times will be inaccurate.
	Results with different $n_{\tt delay}$ will be compared in Section \ref{section:Result}.
\end{remark}

\begin{remark}\rm
	The latent variables are of low dimension so that no truncation in the SVD is performed in our tests,
	i.e. $\epsilon$ is set to zero. In such case, $L=n_{\tt latent}n_{\tt delay}$.
\end{remark}

\subsection{Interpolation for a new parameter in the latent manifold}\label{sec:Interp}
If only one parameter ($d_p=1$) is considered, one can use Lagrangian or linear interpolation,
or manifold interpolation \cite{Zimmermann2021Manifold} to obtain the latent variables at the new parameter value $\widetilde{\bm{\omega}}$
\begin{equation*}
	\bq(t, \widetilde{\bm{\omega}}) = \mathcal{I}(t, \widetilde{\bm{\omega}}) = \mathcal{I}\left(\bq(t, \bm{\omega}_1), \cdots, \bq(t, \bm{\omega}_{n_p})\right).
\end{equation*}
For more parameters, the radial basis function (RBF) interpolation can be employed
\begin{equation*}
	\mathcal{I}(t, \widetilde{\bm{\omega}}) = \sum_{s=1}^{n_p} \bm{w}_s\varphi(\norm{\bm{\omega} - \bm{\omega}_s}),
\end{equation*}
where $\varphi(r) = \varphi(\norm{\bm{\omega} - \bm{\omega}_s})$ is the kernel function,
e.g., the thinplate kernel $\varphi(r) = r^2\ln(r+1)$,
and the weight $\bm{w}_s$ for each parameter $\bm{\omega}_s$ is obtained by solving the least square problem
\begin{equation*}
	\sum_{s=1}^{n_p} \bm{w}_s\varphi(\norm{\bm{\omega}_p - \bm{\omega}_s}) = \bq(t, \bm{\omega}_p),~ p = 1,\cdots, n_p.
\end{equation*}
The RBF interpolation has been studied extensively, and we refer interested readers to \cite{Powell1987Radial} for more details.

\begin{remark}\rm
	The parametric DMD method in this section is similar to the partitioned approach in \cite{Andreuzzi2011A},
	except that we use the autoencoder to generate the latent space rather than the classic POD.
	There are also other parametric DMD approaches,
	such as interpolating the DMD eigenpairs or operators \cite{Huhn2022Parametric},
	using manifold interpolation \cite{Hess2022A}.
	In a very recent work \cite{Conti2022Reduced}, the autoencoder is combined with the SINDy approach for periodic problems.
\end{remark}

%

\subsection{Non-intrusive data-driven RBMs}
Based on the discussions in the last three sections,
we propose the workflow, termed {\tt CAE-PHODMD}, for the construction of the non-intrusive data-driven RBMs
in Algorithms \ref{alg:offline}-\ref{alg:online}, which describe the offline and online stages, respectively.
During the offline stage, we collect suitable snapshots, train the autoencoder, and build the HODMD models for each parameter value.
To predict a full-order solution at a new time and parameter value at the online stage,
we first obtain the latent variables at the new time, then interpolate the new latent variables in the parameter space,
and finally recover the full-order solution by using the decoder.

\begin{algorithm}[H]
	\SetAlgoLined
	\KwIn{Snapshots of the full-order solution at $t_1,\cdots,t_{n_t}$ and $\bm{\omega}_1,\cdots,\bm{\omega}_{n_p}$}
	\KwOut{Autoencoder: $\encop, \decop$; HODMD: $\Lambda, \Theta$ for each parameter value $\bm{\omega}_1,\cdots,\bm{\omega}_{n_p}$}
	Train the autoencoder to get the optimal parameters $\Theta_{\enc}, \Theta_{\dec}$\;
	Compute the latent variables by the encoder 
	$\left[\bq(t_1,\bm{\omega}_1), \dots \bq(t_{n_t},\bm{\omega}_1), \dots\dots, \bq(t_{n_t},\bm{\omega}_{n_p})\right]
	= \left[\encop(\bu_h(t_1,\bm{\omega}_1);\Theta_{\enc}), \dots \encop(\bu_h(t_{n_t},\bm{\omega}_1);\Theta_{\enc}), \dots\dots, \encop(\bu_h(t_{n_t},\bm{\omega}_{n_p});\Theta_{\enc})\right]$\;
	Build the HODMD models for each parameter value $\bm{\omega}_j$ based on the snapshots of the latent variables $\left[\bm{q}(t_1,\bm{\omega}_j), \cdots \bm{q}(t_{n_t},\bm{\omega}_j)\right]$.
	\caption{{\tt CAE-PHODMD}: Offline stage}
	\label{alg:offline}
\end{algorithm}

\begin{algorithm}[H]
	\SetAlgoLined
	\KwIn{A new time $\tilde{t}$ and parameter value $\widetilde{\bm{\omega}}$}
	\KwOut{The full-order solution $\bu(\tilde{t}, \widetilde{\bm{\omega}})$}
	Predict in time using the HODMD: $\mathcal{Q}(\bm{\omega}_1),\cdots,\mathcal{Q}(\bm{\omega}_{n_p}) \Rightarrow \left[\bm{q}(\tilde t,\bm{\omega}_1), \cdots, \bm{q}(\tilde t,\bm{\omega}_{n_p}) \right]$\;
	Interpolate for a new parameter value: $\left[\bm{q}(\tilde t,\bm{\omega}_1), \cdots, \bm{q}(\tilde t,\bm{\omega}_{n_p}) \right] \Rightarrow \bm{q}(\tilde t,\widetilde{\bm{\omega}})$\;
	Reconstruct by the decoder: $\bu_h(\tilde{t}, \widetilde{\bm{\omega}}) = \decop(\bm{q}(\tilde t,\widetilde{\bm{\omega}});\Theta_{\dec})$.
	\caption{{\tt CAE-PHODMD}: Online stage}
	\label{alg:online}
\end{algorithm}


\section{Numerical results}\label{section:Result}
This section presents numerical tests on typical parametrized PDEs.
The following relative errors at a given testing time and parameter value $(t_{\tt test}, \bm{\omega}_{\tt test})$ will be evaluated
\begin{align*}
	&\epsilon_{\tt CAE} (t_{\tt test}, \bm{\omega}_{\tt test})=
	\dfrac{{\norm{\bu_h(t_{\tt test}, \bm{\omega}_{\tt test})
				- \decop(\encop(\bu_h(t_{\tt test}, \bm{\omega}_{\tt test}); \Theta_{\enc}); \Theta_{\dec})}_2}}
		{{\norm{\bu_h(t_{\tt test}, \bm{\omega}_{\tt test})}_2}}, \\
	&\epsilon_{\tt latent} (t_{\tt test}, \bm{\omega}_{\tt test})=
	\dfrac{{\norm{\bq(t_{\tt test}, \bm{\omega}_{\tt test})
				- \mathcal{I}\left(t_{\tt test}, \bm{\omega}_{\tt test}\right)}_2}}
		{{\norm{\bq(t_{\tt test}, \bm{\omega}_{\tt test})}_2}}, \\
	&\epsilon_{\tt CAE-PHODMD} (t_{\tt test}, \bm{\omega}_{\tt test})=
	\dfrac{{\norm{\bu_h(t_{\tt test}, \bm{\omega}_{\tt test})
				- \decop(\mathcal{I}\left(t_{\tt test}, \bm{\omega}_{\tt test}\right); \Theta_{\dec})}_2}}
	{{\norm{\bu_h(t_{\tt test}, \bm{\omega}_{\tt test})}_2}},
\end{align*}
where $\bq(t_{\tt test}, \bm{\omega}_{\tt test}) = \encop(\bu_h(t_{\tt test}, \bm{\omega}_{\tt test}); \Theta_{\enc})$
is the vector of latent variables obtained by the encoder,
and $\mathcal{I}\left(t_{\tt test}, \bm{\omega}_{\tt test}\right)
=\mathcal{I}\left(\bq(t_{\tt test}, \bm{\omega}_{1}), \cdots, \bq(t_{\tt test}, \bm{\omega}_{n_p})\right)$
is the interpolation of the latent variables at the new parameter value.
The errors of the whole testing set are defined as 
\begin{align*}
	&E_{\tt CAE} = \dfrac{1}{N_{t_{\tt test}}N_{\bm{\omega}_{\tt test}}}
	\sum\limits_{t_{\tt test}}\sum\limits_{\bm{\omega}_{\tt test}}
	\epsilon_{\tt CAE} (t_{\tt test}, \bm{\omega}_{\tt test}), \\
	&E_{\tt latent} = \dfrac{1}{N_{t_{\tt test}}N_{\bm{\omega}_{\tt test}}}
	\sum\limits_{t_{\tt test}}\sum\limits_{\bm{\omega}_{\tt test}}
	\epsilon_{\tt latent} (t_{\tt test}, \bm{\omega}_{\tt test}), \\
	&E_{\tt CAE-PHODMD} = \dfrac{1}{N_{t_{\tt test}}N_{\bm{\omega}_{\tt test}}}
	\sum\limits_{t_{\tt test}}\sum\limits_{\bm{\omega}_{\tt test}}
	\epsilon_{\tt CAE-PHODMD} (t_{\tt test}, \bm{\omega}_{\tt test}),
\end{align*}
where $N_{t_{\tt test}}$ and $N_{\bm{\omega}_{\tt test}}$ are the numbers of the testing times and parameter values, respectively.
Our implementation of the autoencoder is based on PyTorch library \cite{PyTorch2019}.
In all the tests, the mini-batch ADAM optimizer with an initial learning rate $0.001$ and batch size $32$ is adopted for training.
And the StepLR scheduler with step size $50$ and decay rate $0.95$ is used,
so that the learning rate is $\lambda=0.001\times 0.95^{\lfloor n/50\rfloor}$, where $n$ is the current number of epochs.
In the training of the CAE, early stopping technique is used to avoid overfitting.
To be specific, the training is terminated if the best reconstruction error on the validation set has not been improved for $100$ epochs.
The linear interpolation provided by SciPy package is used in all the tests.
There is no difference if using higher-order interpolation.

\subsection{1D Burgers' equation}\label{sec:1DBurgers}
\subsubsection{Setup}
This example is used to verify our approach on the 1D Burgers' equation
\begin{equation*}
	u_t + \left(\frac12u^2\right)_x = \frac{1}{Re} u_{xx},~ x\in[0,2],~ u(0,t)=u(2,t)=0,
\end{equation*}
with the exact solution \cite{Maulik2021Reduced}
\begin{equation*}
	u(x,t;Re) = \dfrac{\dfrac{x}{t+1}}{1 + \sqrt{\dfrac{t+1}{t_0} \exp(Re \dfrac{x^2}{4t+4})}},\quad t_0=\exp(Re/8),
\end{equation*}
where $Re$ is the Reynolds number.
In this test, the exact solutions with $128$ uniform spatial degrees of freedom serve as the full-order solutions.
The datasets consist of three parts:
the training set comprises $10$ $Re$ uniform in $[100,800]$ with $101$ times uniform in $[0,2]$,
the validation set comprises $4$ random $Re$ with $20$ random times in the same parameter and time domains,
and the testing set comprises $2$ random $Re$ with $10$ random times in the same parameter and time domains.
The kernel size in all the convolutional and deconvolutional layers is $5$ with stride and padding as $2$.

\begin{table}[hbt!]
	\centering
	\begin{tabular}{c|r|r}
			\toprule
			\multicolumn{3}{c}{encoder} \\ \hline
		   layer & input shape & output shape \\ \hline
		   {\tt Conv1d with SiLU} & $1\times 128$ & $32\times 64$ \\
		   {\tt Conv1d with SiLU} & $32\times 64$ & $32\times 32$ \\
		   {\tt Conv1d with SiLU} & $32\times 32$ & $32\times 16$ \\
		   {\tt Conv1d with SiLU} & $32\times 16$ & $32\times 8$ \\
		   {\tt Conv1d with SiLU} & $32\times 8$ & $32\times 4$ \\
		   {\tt Conv1d with SiLU} & $32\times 4$ & $32\times 2$ \\
		   {\tt Flatten} & $32\times 2$ & $ 64 $ \\
		   {\tt Linear with SiLU} & $64$ & $11$ \\
		   {\tt Linear with SiLU} & $11$ & $2$ \\
			\toprule
			\multicolumn{3}{c}{decoder} \\ \hline
			layer & input shape & output shape \\ \hline		
			{\tt Linear with SiLU} & $2$ & $11$ \\
			{\tt Linear with SiLU} & $ 11 $ & $64$ \\
			{\tt Unflatten} & $64$ & $32\times 2$ \\
			{\tt ConvTranspose1d with SiLU} & $32\times 2$ & $32\times 4$ \\
			{\tt ConvTranspose1d with SiLU} & $32\times 4$ & $32\times 8$ \\
			{\tt ConvTranspose1d with SiLU} & $32\times 8$ & $32\times 16$ \\
			{\tt ConvTranspose1d with SiLU} & $32\times 16$ & $32\times 32$ \\
			{\tt ConvTranspose1d with SiLU} & $32\times 32$ & $32\times 64$ \\
			{\tt ConvTranspose1d with SiLU} & $32\times 64$ & $1\times 128$ \\
			\bottomrule
		\end{tabular}
	\caption{1D Burgers' equation: The architecture of the CAE with the best reconstruction error on the validation set during the grid search.
	The weight decay is $10^{-11}$ in the training.}
	\label{tab:1DBurgers_arch}
\end{table}

\subsubsection{Results}
The grid search is performed to find the best architecture of the CAE,
with different weight decay $10^{-8}$, $10^{-9}$, $10^{-10}$, $10^{-11}$,
number of convolutional layers $4,5,6$,
and the dimension of the latent space $n_{\tt latent}=2$, $4$, $6$, $8$, $10$.
The architecture that obtains the best reconstruction error on the validation set is chosen as the final CAE model, shown in Table \ref{tab:1DBurgers_arch}.
The CAE models with other $n_{\tt latent}$ are also used for comparison.

\begin{figure}[hbt!]
	\centering
	\includegraphics[width=0.6\textwidth]{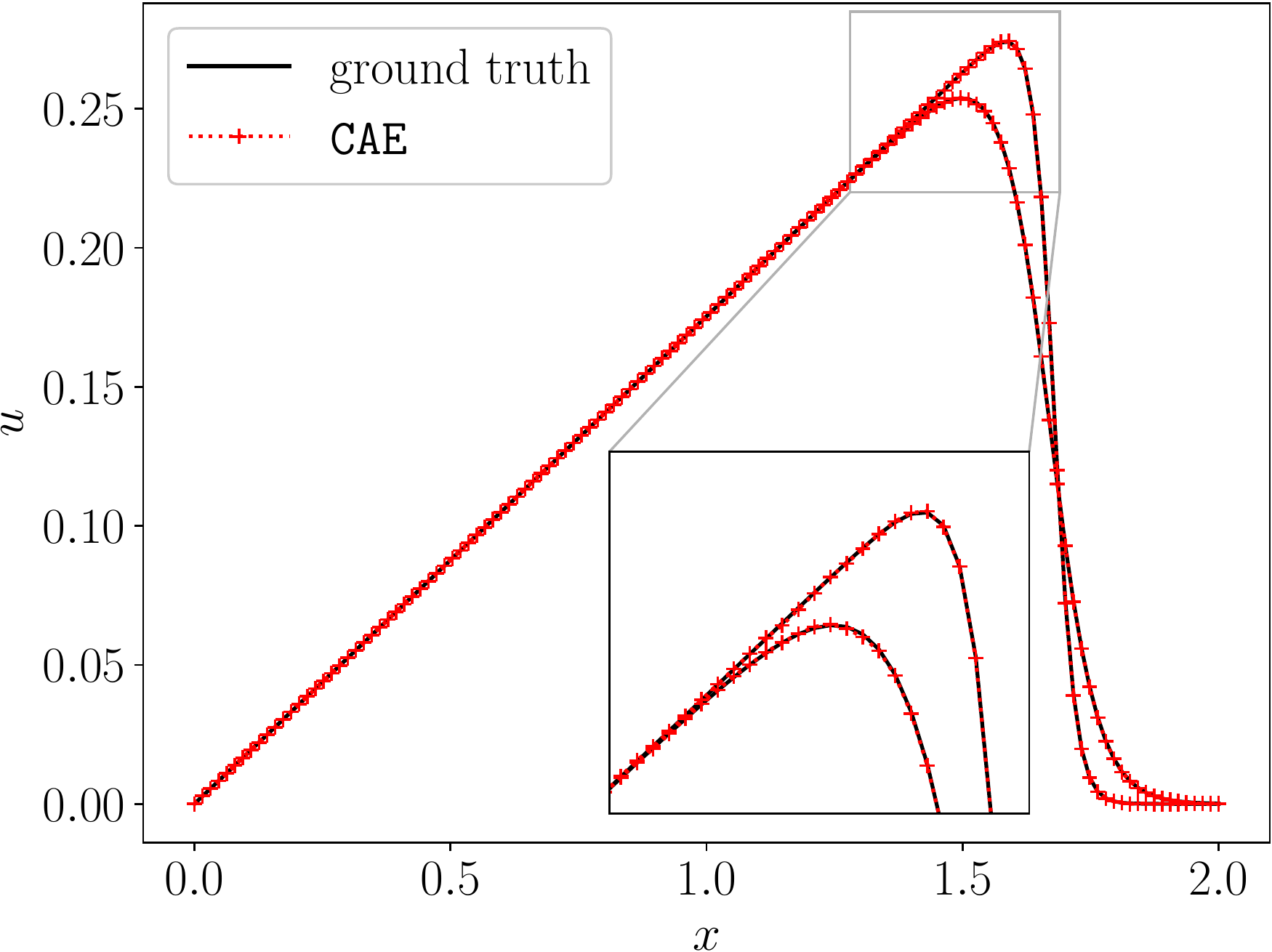}
	\caption{1D Burgers' equation.
	The reconstructed full-order solutions in the testing set for $Re=290.8$, $646.0$ at $t=1.85$ by the CAE with $n_{\tt latent} = 2$.}
	\label{fig:1DBurgers_CAE}
\end{figure}

Figure \ref{fig:1DBurgers_CAE} shows two specific reconstructed full-order solutions in the testing set for $Re=290.8$, $646.0$ at $t=1.85$,
and the corresponding errors are $\epsilon_{\tt CAE} = 1.002\times 10^{-3}$, $7.679\times 10^{-4}$.
It is seen that $n_{\tt latent} = 2$ suffices to obtain accurate reconstructions without oscillations,
which indicates the high efficiency of the CAE.

\begin{figure}[hbt!]
    \centering
    \includegraphics[width=0.8\textwidth]{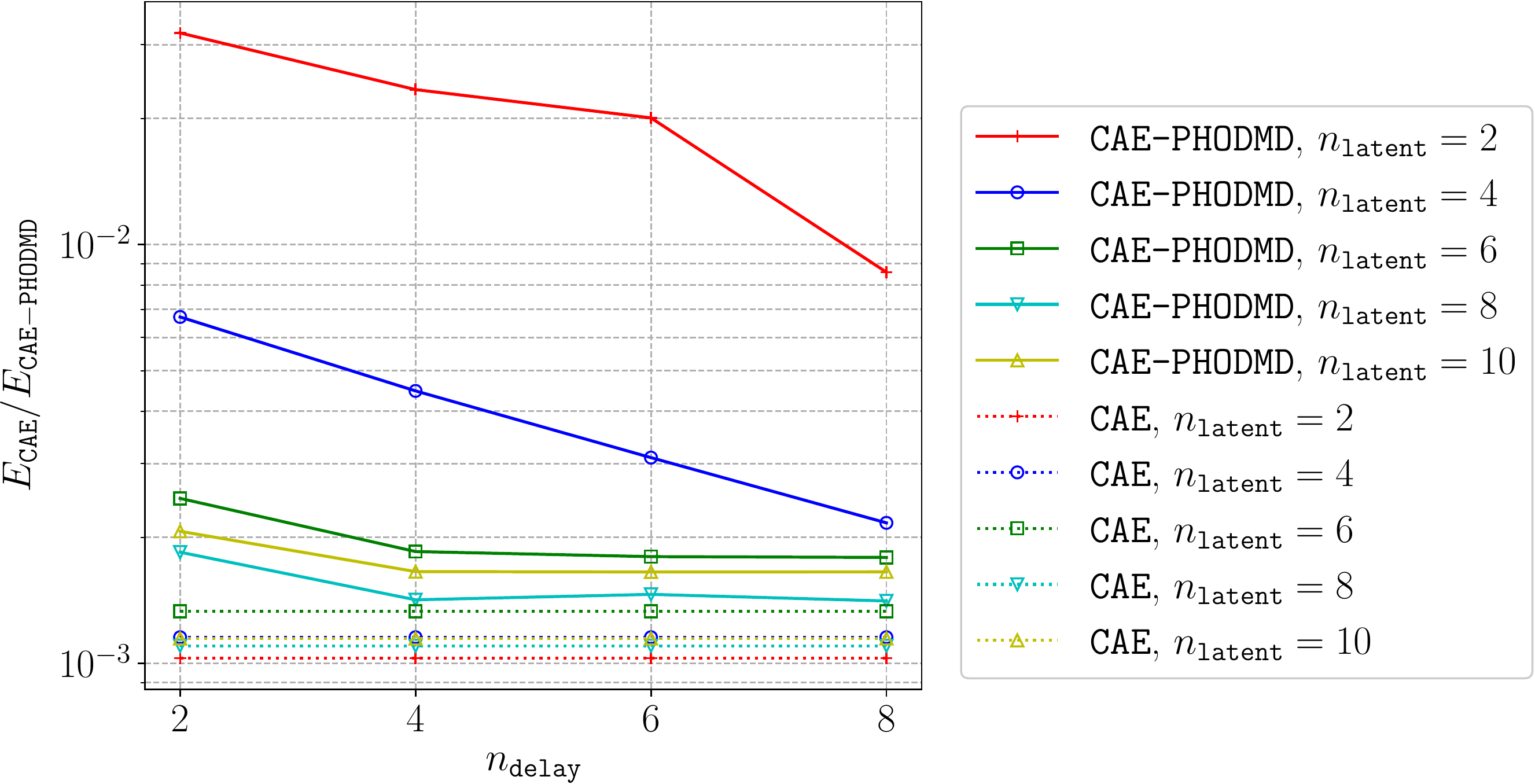}
    \caption{1D Burgers' equation.
    The errors of the reconstructed full-order solution in the testing set by the CAE $E_{\tt CAE}$ with different $n_{\tt latent}$, 
    and prediction errors by {\tt CAE-PHODMD}, $E_{\tt CAE-PHODMD}$ w.r.t. the number of time-delay embedding $n_{\tt delay}$.}
    \label{fig:1DBurgers_err_CAE_PHODMD}
\end{figure}

The performance of {\tt CAE-PHODMD} depends on the dimension of the latent space $n_{\tt latent}$,
and also the number of time-delay embedding $n_{\tt delay}$ used in the HODMD.
Figure \ref{fig:1DBurgers_err_CAE_PHODMD} presents the errors $E_{\tt CAE-PHODMD}$ with different $n_{\tt latent}, n_{\tt delay}$,
and the errors $E_{\tt CAE}$ with different $n_{\tt latent}$.
For $n_{\tt latent} = 2, 4$, the errors of {\tt CAE-PHODMD} decay as $n_{\tt delay}$ increases,
while for $n_{\tt latent} = 6, 8, 10$, the errors remain constant for $n_{\tt delay} \geqslant 4$,
because there are sufficient time-delay embeddings to obtain accurate results.
One can see that $E_{\tt CAE-PHODMD}$ is generally larger than $E_{\tt CAE}$,
since the HODMD and interpolation introduce extra errors.

\begin{figure}[hbt!]
	\centering
	\begin{subfigure}[b]{0.49\textwidth}
		\centering
		\includegraphics[width=1.0\textwidth]{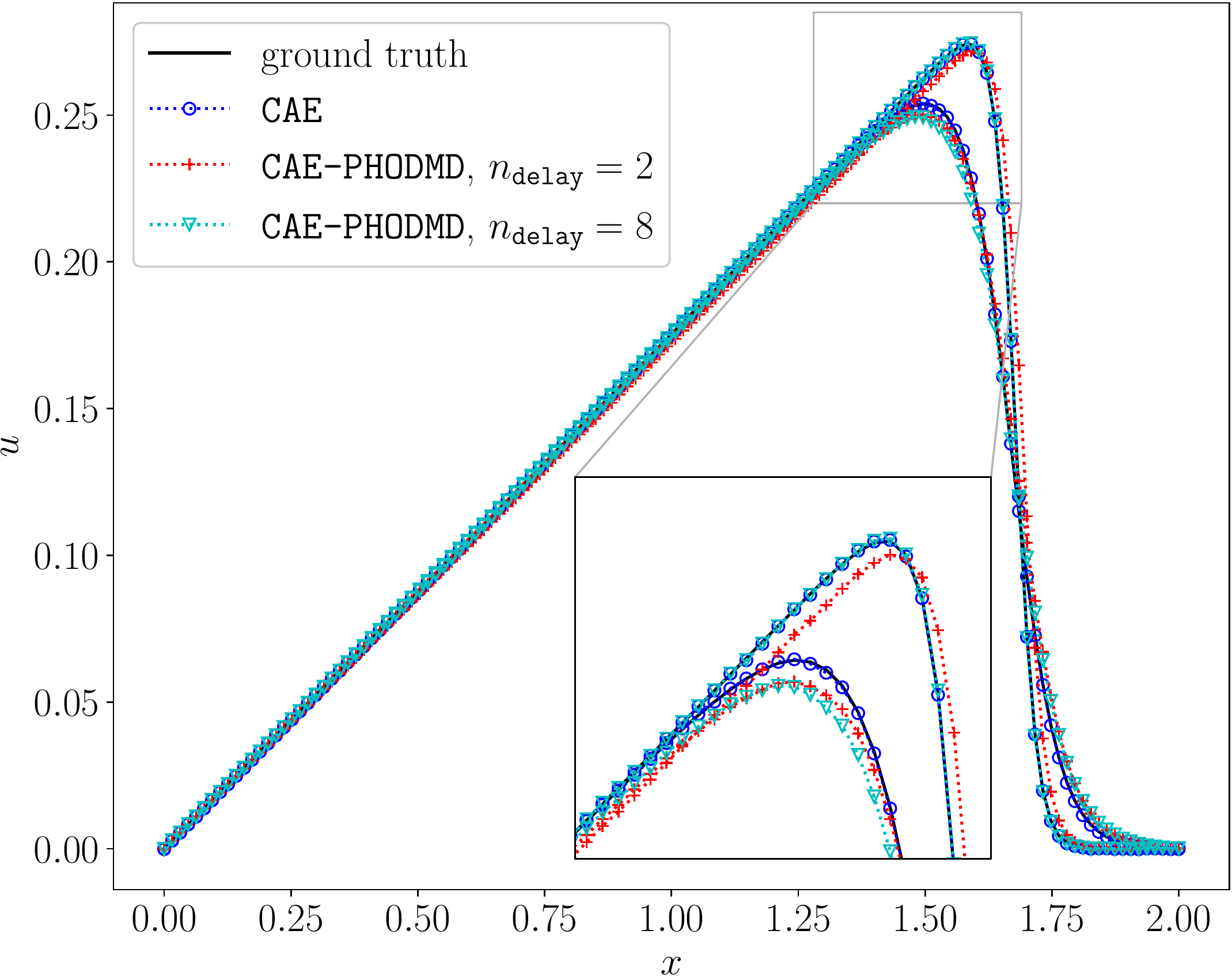}
		\caption{$n_{\tt latent} = 2$}
	\end{subfigure}
	\begin{subfigure}[b]{0.49\textwidth}
		\centering
		\includegraphics[width=1.0\textwidth]{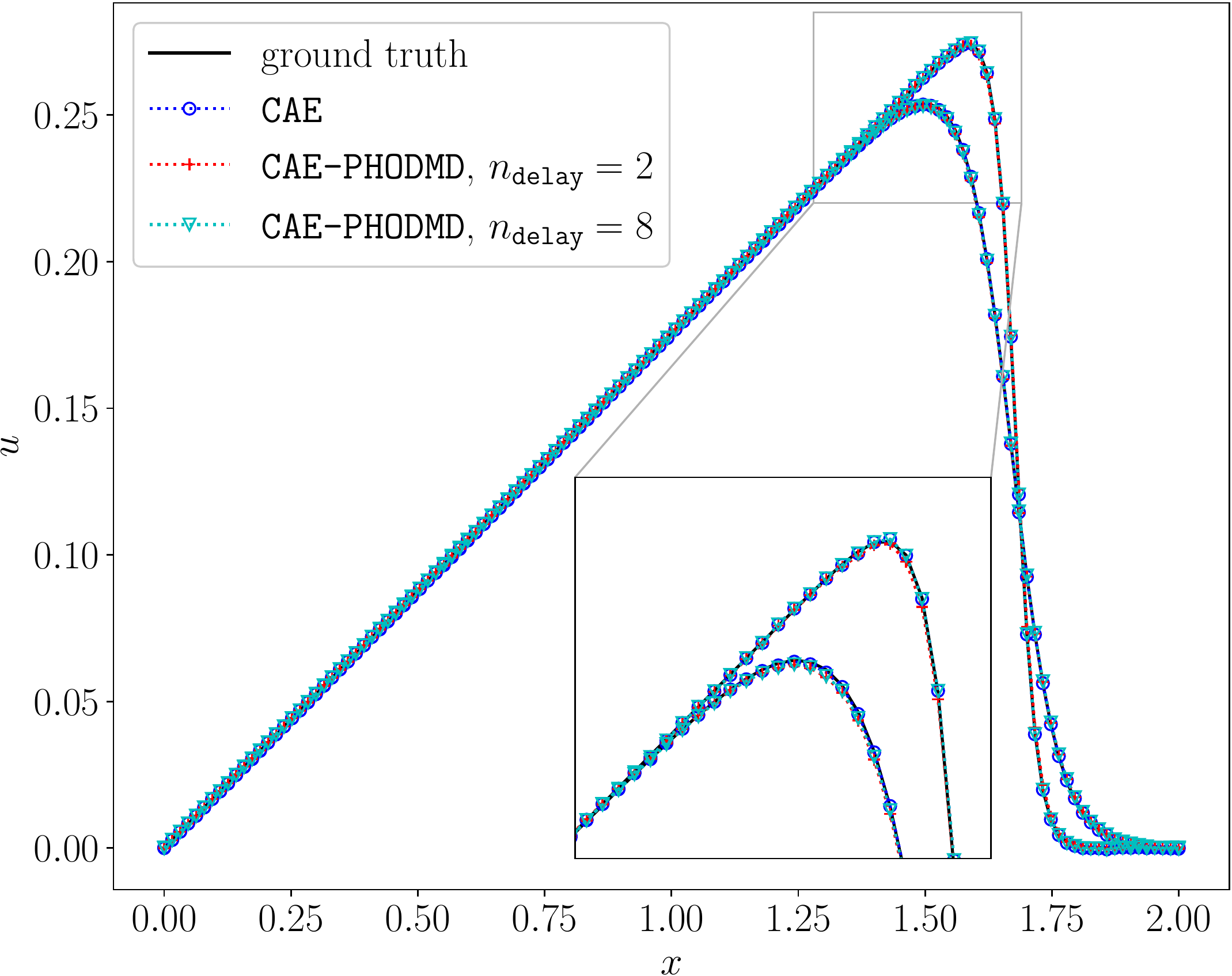}
		\caption{$n_{\tt latent} = 6$}
	\end{subfigure}
	\caption{1D Burgers' equation. The reconstructed full-order solution for $Re=290.8$, $646.0$ at $t=1.85$ by the CAE and {\tt CAE-PHODMD}.}
	\label{fig:1DBurgers_CAE_PHODMD_u_latent2}
\end{figure}

The predicted full-order solutions for $Re=290.8$, $646.0$ at $t=1.85$
obtained by {\tt CAE-PHODMD} with $n_{\tt latent} = 2, 6$ and $n_{\tt delay} = 2, 8$
are plotted in Figure \ref{fig:1DBurgers_CAE_PHODMD_u_latent2}.
The left figure shows that increasing $n_{\tt delay}$ improves the accuracy for $n_{\tt latent}=2$,
and from the right figure one observes that
$n_{\tt delay}=2$ suffices to make the results close to the reconstructed full-order solutions obtained by the CAE.
In this case, $n_{\tt latent}=6$ with $n_{\tt delay}=2$ is enough to ensure accurate prediction in the testing set.

\begin{figure}[hbt!]
	\centering
	\includegraphics[width=0.6\textwidth]{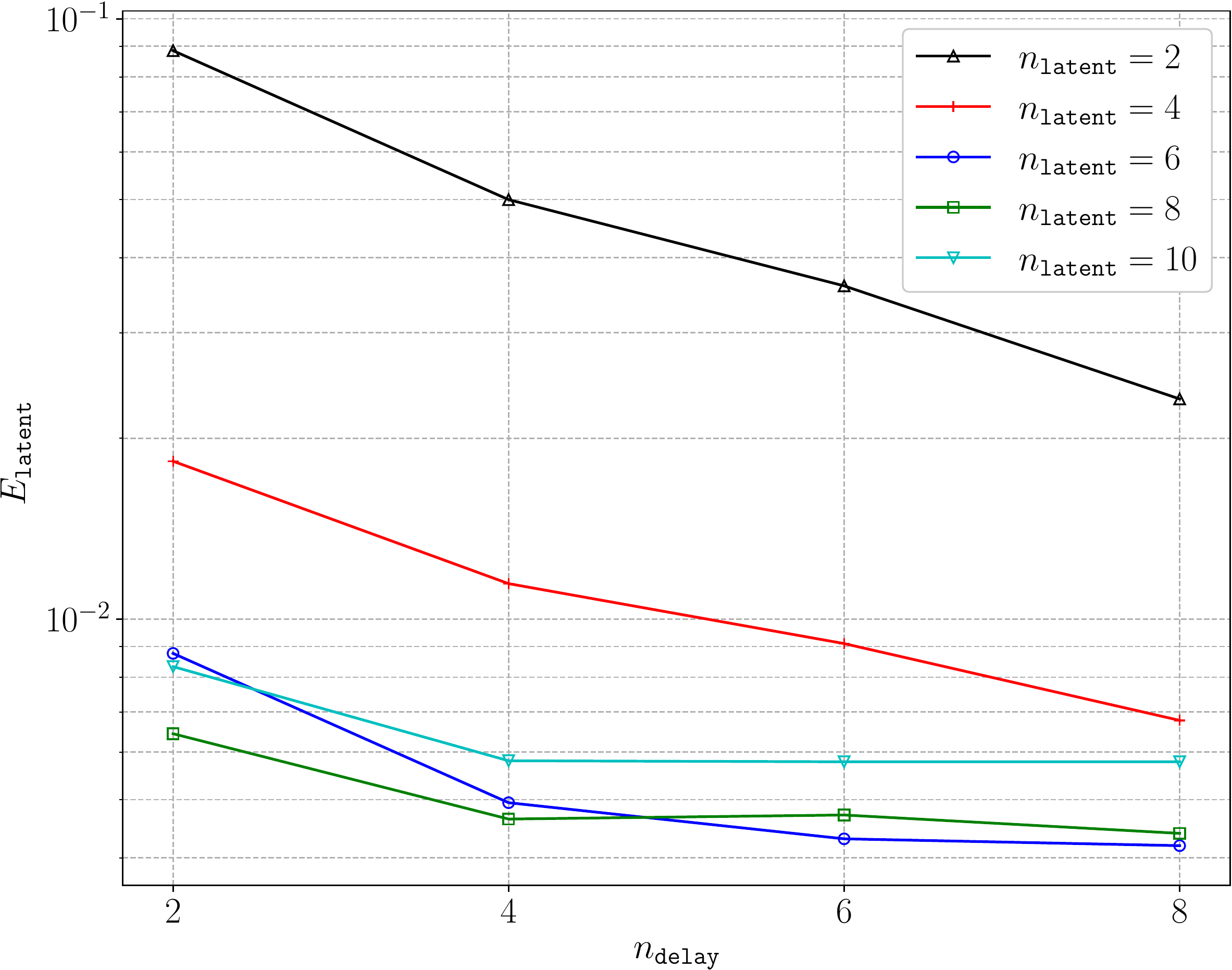}
	\caption{1D Burgers' equation. The errors in the latent variables in the testing set by {\tt CAE-PHODMD},
		$E_{\tt latent}$ w.r.t. the number of time-delay embedding $n_{\tt delay}$.}
	\label{fig:1DBurgers_err_latent}
\end{figure}

\begin{figure}[hbt!]
	\centering
	\begin{subfigure}[b]{0.49\textwidth}
		\centering
		\includegraphics[width=1.0\textwidth]{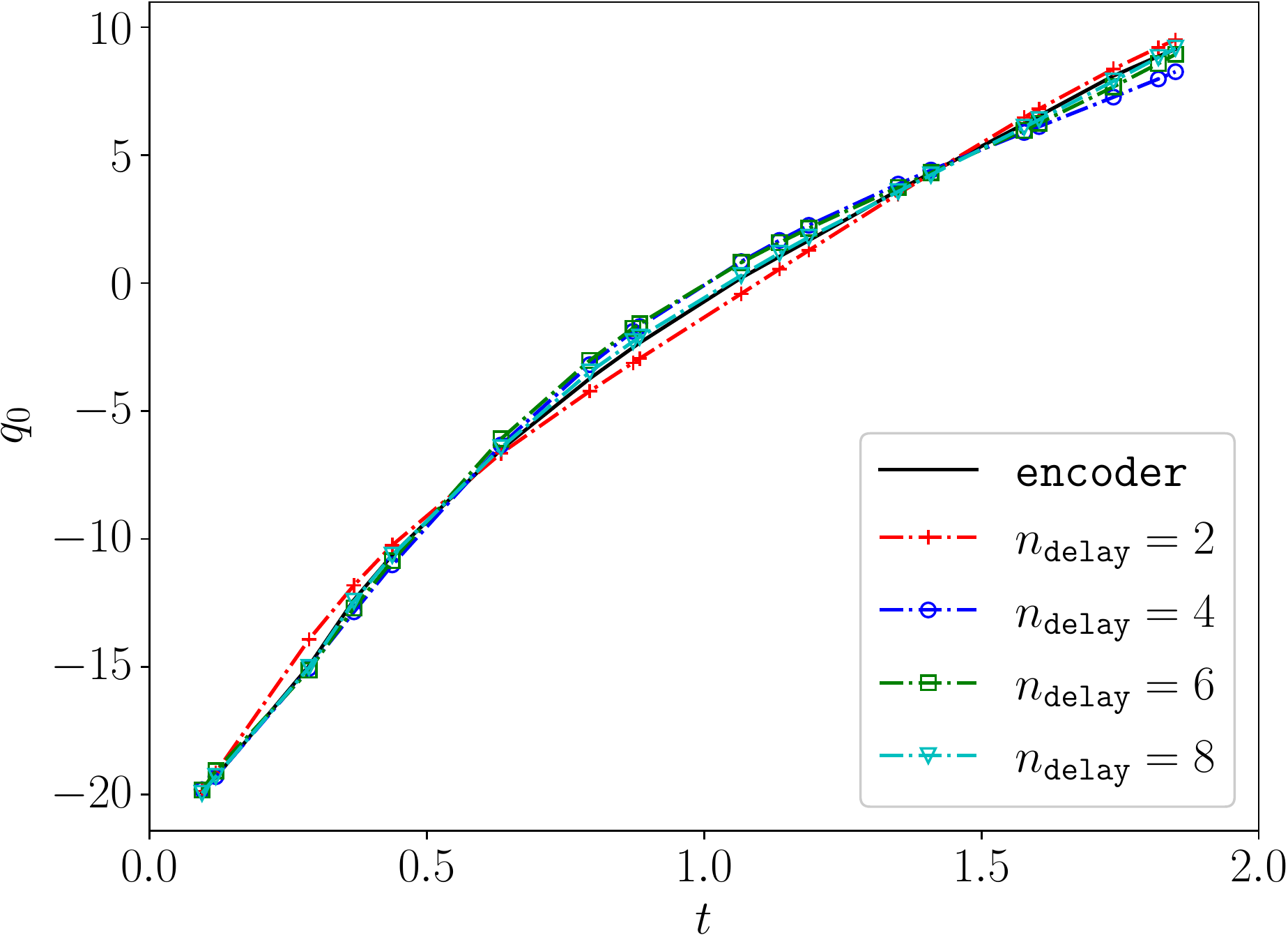}
		\caption{The first latent variable.}
	\end{subfigure}
	\begin{subfigure}[b]{0.49\textwidth}
		\centering
		\includegraphics[width=1.0\textwidth]{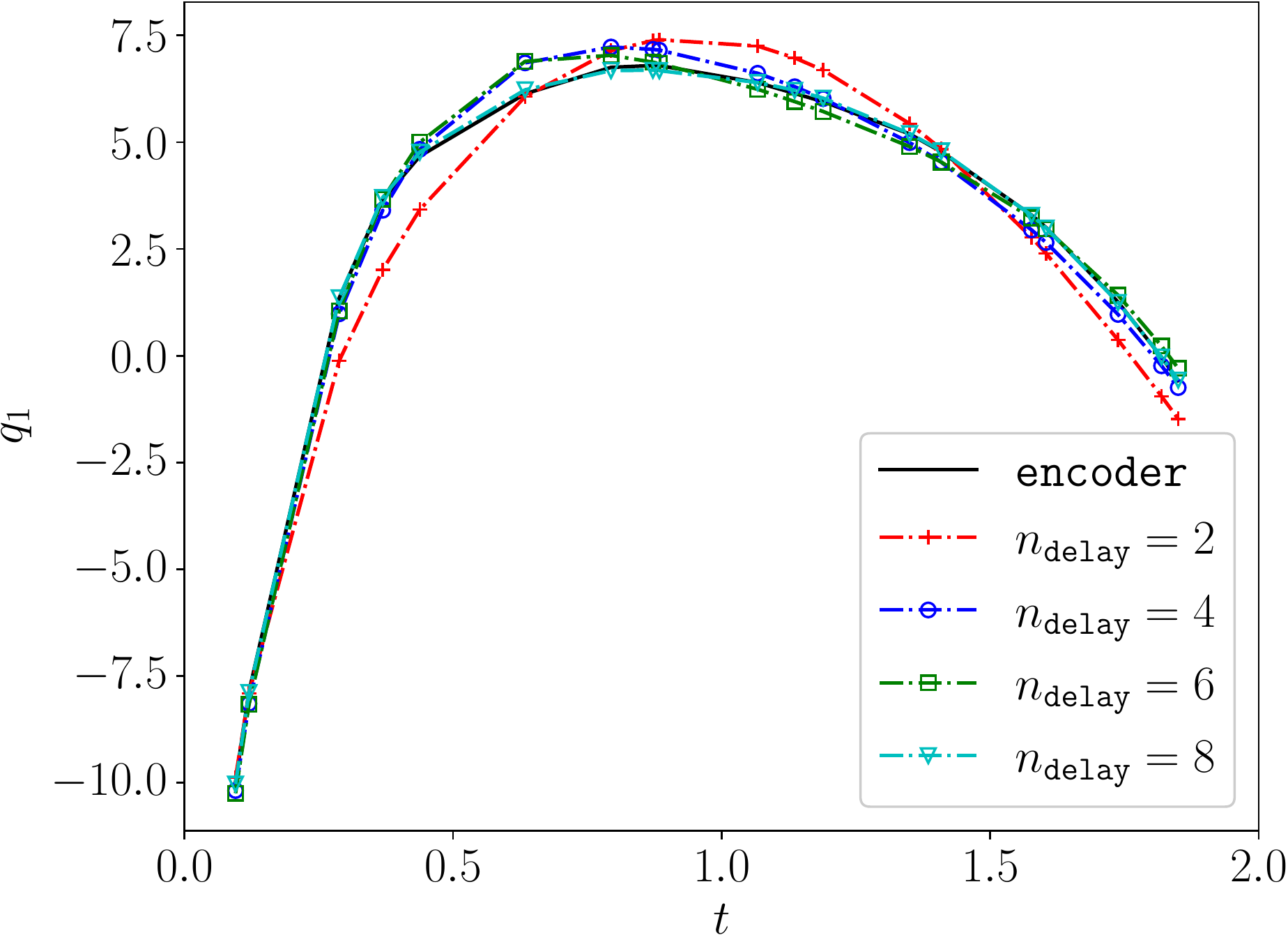}
		\caption{The second latent variable.}
	\end{subfigure}
	\caption{1D Burgers' equation.
		The temporal evolution of the latent variables for $Re=646.0$ obtained by {\tt CAE-PHODMD}.}
	\label{fig:1DBurgers_latent}
\end{figure}

To investigate the errors in the latent space,
Figure \ref{fig:1DBurgers_err_latent} shows the errors $E_{\tt latent}$ obtained by {\tt CAE-PHODMD}.
Similar to Figure \ref{fig:1DBurgers_err_CAE_PHODMD},
the errors decrease as $n_{\tt delay}$ increases for $n_{\tt latent} = 2,4$,
and do not improve for $n_{\tt latent} = 6,8,10$ when $n_{\tt delay} \geqslant 4$,
thus the error of the prediction by {\tt CAE-PHODMD} is dominated by the error in the latent space in this case.
The temporal evolution of the latent variables with $n_{\tt latent}=2$ for $Re=646.0$ obtained by {\tt CAE-PHODMD} is given in Figure \ref{fig:1DBurgers_latent}.
It also shows that the HODMD becomes more accurate as the number of time-delay embedding $n_{\tt delay}$ increases.

\subsection{2D Rayleigh-B\'enard convection}\label{sec:2DRBC}
\subsubsection{Setup}
This test will verify the effectiveness and performance of our approach on the 2D Rayleigh-B\'enard convection problem.
The full-order solutions are obtained by solving the 2D incompressible Navier-Stokes equations
\begin{equation*}
	\left\{~
	\begin{aligned}
		&\nabla\cdot\bv=0, \\
		&\bv_t + \bv\cdot \nabla \bv = -\nabla p + \sqrt{\frac{Pr}{Ra}} \Delta \bv + T \bm{e}_y, \\
		&T_t + \bv\cdot \nabla T = \frac{1}{\sqrt{Pr Ra}} \Delta T,
	\end{aligned}
	 \right.
\end{equation*}
where $\bv, T, p$ are the dimensionless velocity, temperature and pressure, respectively.
The Rayleigh number $Ra$ and the Prandtl number $Pr$ are two dimensionless parameters controlling the flow,
and $Pr$ is fixed as $0.71$ in this test.
The computational domain is $[0,2]\times[0,1]$ with the wall boundary conditions for the velocity at all boundaries.
The high and low temperatures are specified at the lower and upper boundaries, respectively, $T(y=0) = 1$, $T(y=1) = 0$,
and the zero Neumann boundary conditions are used for the temperature at the left and right boundaries.
The solution vector $\bu=(\bv, T)$ consists of the velocity and temperature,
and semi-implicit scheme with Taylor-Hood element on the $128\times 64$ uniform quadrilateral mesh
is employed to compute the full-order solutions.
In all the simulations, the initial data are chosen as the steady state of $Ra=10^{4}$.
The kernel size used in all the convolutional and deconvolutional layers is $5$ with the stride and padding as $2$.

The datasets used in this test consist of three parts.
The training set comprises the snapshots at $Ra=10^6$, $2\times 10^6$, $4\times 10^6$, $5\times 10^6$, $7\times 10^6$, $9\times 10^6$,
and $200$ times uniform in $[0,40]$.
The validation set comprises the snapshots at $Ra=3\times 10^6$, $8\times 10^6$,
and $20$ random times in $[0,40]$.
The testing set comprises the snapshots at $Ra=6\times 10^6$,
and $10$ random times in $[0,48]$.

\begin{table}[hbt!]
	\centering
	\begin{tabular}{c|r|r}
		\toprule
		\multicolumn{3}{c}{encoder} \\ \hline
		layer & input shape & output shape \\ \hline
		{\tt Conv2d with SiLU} & $3\times 64\times 128$ & $16\times 32\times 64$ \\
		{\tt Conv2d with SiLU} & $16\times 32\times 64$ & $16\times 16\times 32$ \\
		{\tt Conv2d with SiLU} & $16\times 16\times 32$ & $16\times 8\times 16$ \\
		{\tt Conv2d with SiLU} & $16\times 8\times 16$ & $16\times 4\times 8$ \\
		{\tt Flatten} & $16\times 4\times 8$ & $ 512 $ \\
		{\tt Linear with SiLU} & $ 512 $ & $71$ \\
		{\tt Linear with SiLU} & $ 71 $ & $10$ \\
		\toprule
		\multicolumn{3}{c}{decoder} \\ \hline
		layer & input shape & output shape \\ \hline		
		{\tt Linear with SiLU} & $10$ & $71$ \\
		{\tt Linear with SiLU} & $ 71 $ & $512$ \\
		{\tt Unflatten} & $ 512 $ & $16\times 4\times 8$ \\
		{\tt ConvTranspose2d with SiLU} & $16\times 4\times 8$ & $16\times 8\times 16$ \\
		{\tt ConvTranspose2d with SiLU} & $16\times 8\times 16$ & $16\times 16\times 32$ \\
		{\tt ConvTranspose2d with SiLU} & $16\times 16\times 32$ & $16\times 32\times 64$ \\
		{\tt ConvTranspose2d with SiLU} & $16\times 32\times 64$ & $3\times 64\times 128$ \\
		\bottomrule
	\end{tabular}
	\caption{2D Rayleigh-B\'enard convection: The architecture of the CAE during the grid search with the best reconstruction error on the validation set.
	The weight decay is $10^{-11}$ in the training.}
	\label{tab:2DRBC_arch}
\end{table}

\subsubsection{Results}
To find a preferred architecture, we perform a grid search with different weight decay $10^{-8}$, $10^{-9}$, $10^{-10}$, $10^{-11}$,
number of convolutional layers $4$, $5$, $6$, and the dimension of the latent space $n_{\tt latent}=4$, $6$, $8$, $10$.
The architecture with the best reconstruction error on the validation set is chosen as the final CAE model, shown in Table \ref{tab:2DRBC_arch}.

\begin{figure}[hbt!]
	\centering
	\includegraphics[width=1.0\textwidth]{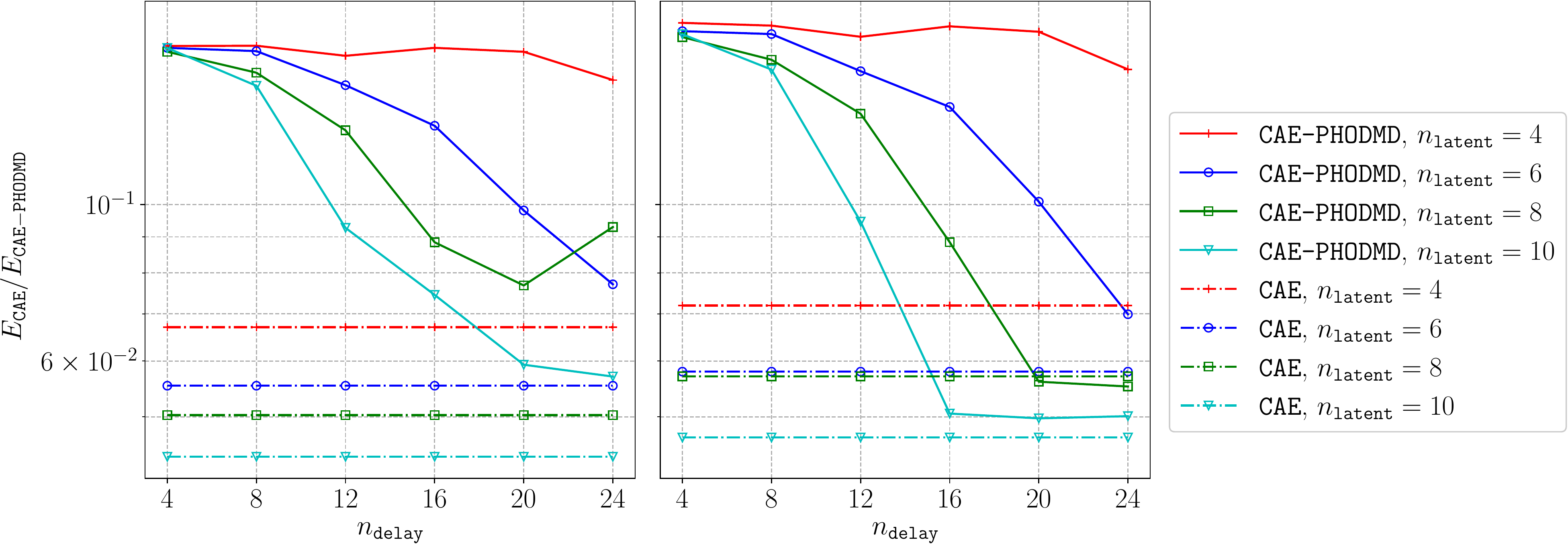}
	\caption{2D Rayleigh-B\'enard convection.
	The reconstruction errors $E_{\tt CAE}$ and prediction errors $E_{\tt CAE-PHODMD}$ w.r.t. the number of the time-delay embedding $n_{\tt delay}$ for different dimensions of the latent space $n_{\tt latent}$.
	The left figure is for all testing times $t\in[0,48]$,
	while the right only for interpolation in time $t\in[0,40]$.
}
	\label{fig:2DRBC_err_u_CAE_PHODMD}
\end{figure}

Figure \ref{fig:2DRBC_err_u_CAE_PHODMD} presents the errors $E_{\tt CAE}$ and $E_{\tt CAE-PHODMD}$ with different $n_{\tt latent}$ and $n_{\tt delay}$.
The left and right figures show the errors for all the testing times and interpolation in time, i.e., $t\in[0.48]$ and $t\in[0.40]$, respectively.
It is observed that the errors in the left figure are generally larger than in the right,
and the errors $E_{\tt CAE-PHODMD}$ are larger than $E_{\tt CAE}$.
One can see that nearly all the errors become smaller when using larger $n_{\tt latent}$,
and $E_{\tt CAE-PHODMD}$ decreases as $n_{\tt delay}$ increases,
then remains constant for $n_{\tt latent}=8,10$ when $n_{\tt delay}\geqslant 20$.
In the left figure, the error $E_{\tt CAE-PHODMD}$ with $n_{\tt latent}=8$ increases at $n_{\tt delay}=24$,
which is due to the extrapolation error.
The one on the right keeps constant when the samples in $t\in[40,48]$ are excluded.

%
%

\begin{figure}[hbt!]
	\centering
	\begin{subfigure}[b]{0.45\textwidth}
		\centering
		\includegraphics[width=1.0\textwidth]{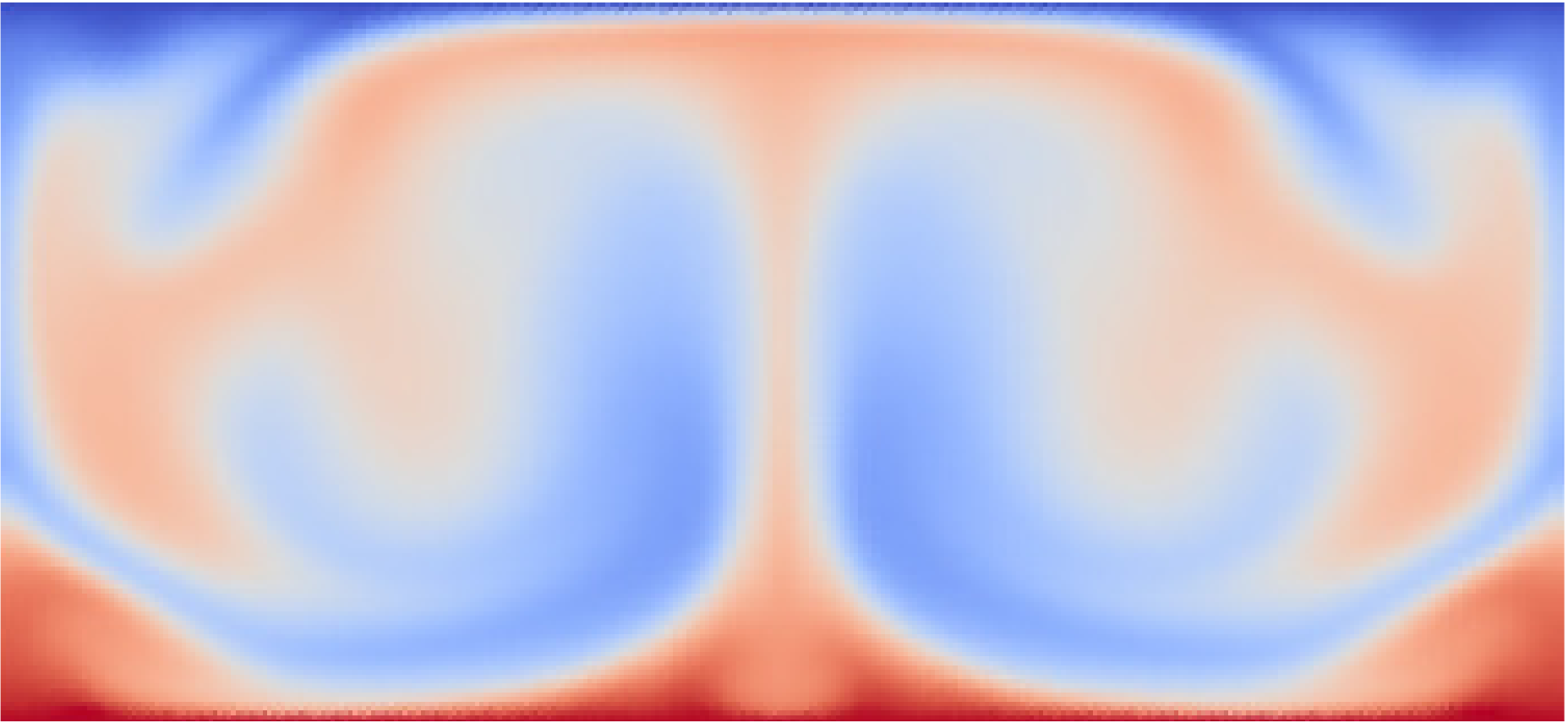}
		\caption{ground truth}
	\end{subfigure}
	\begin{subfigure}[b]{0.45\textwidth}
		\centering
		\includegraphics[width=1.0\textwidth]{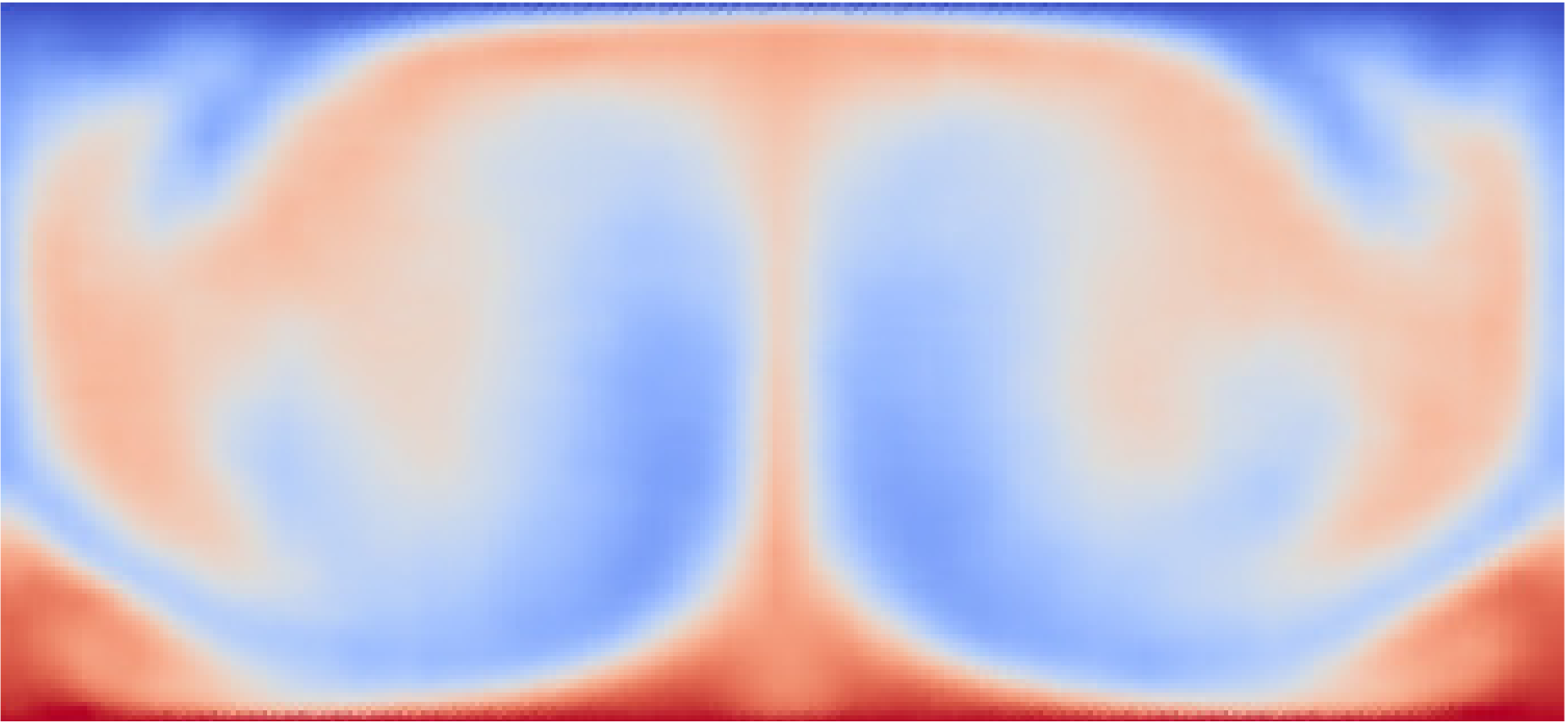}
		\caption{$\epsilon_{\tt CAE}=3.779\times 10^{-2}$}
	\end{subfigure}
	
	\begin{subfigure}[b]{0.45\textwidth}
		\centering
		\includegraphics[width=1.0\textwidth]{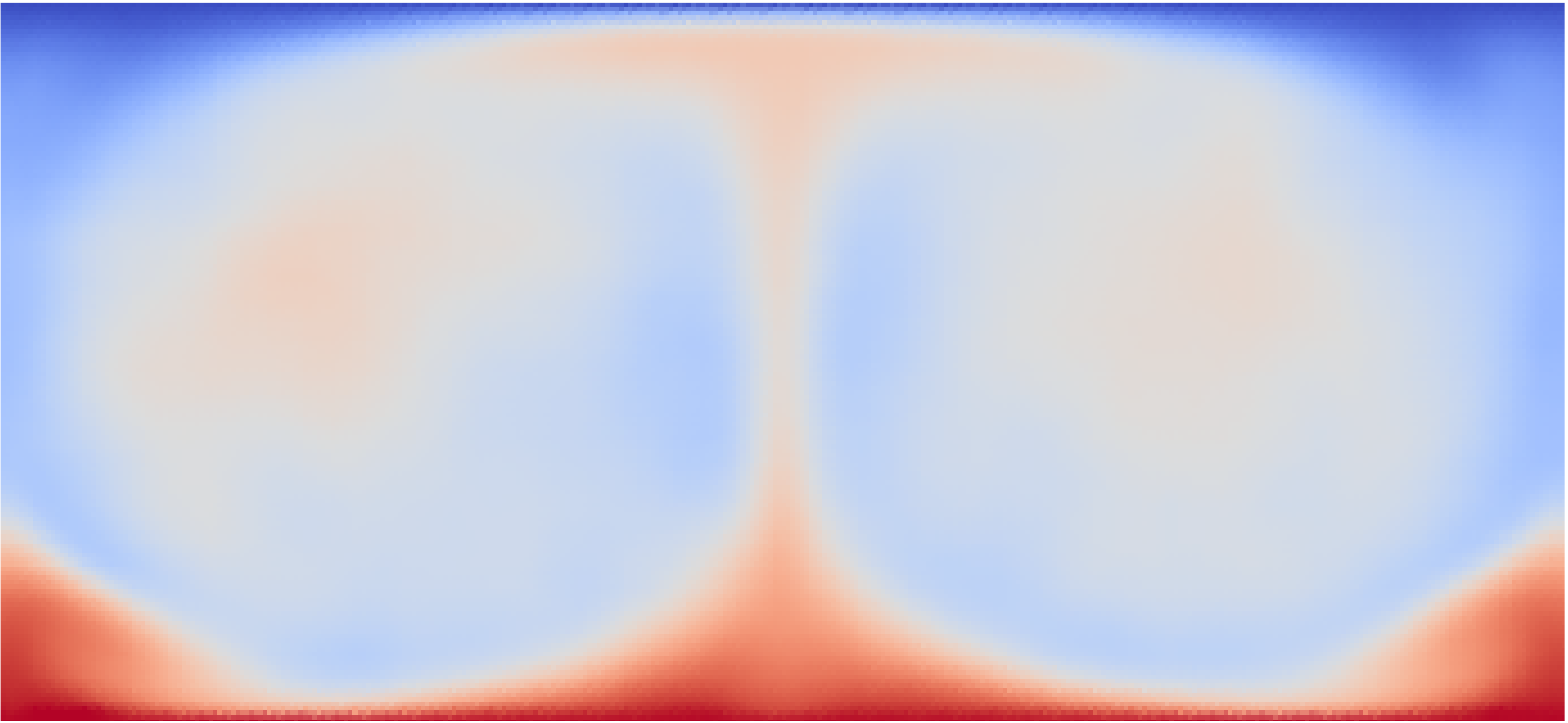}
		\caption{$\epsilon_{\tt CAE-PHODMD}=1.930\times 10^{-1}, n_{\tt delay}=4$}
	\end{subfigure}
	\begin{subfigure}[b]{0.45\textwidth}
		\centering
		\includegraphics[width=1.0\textwidth]{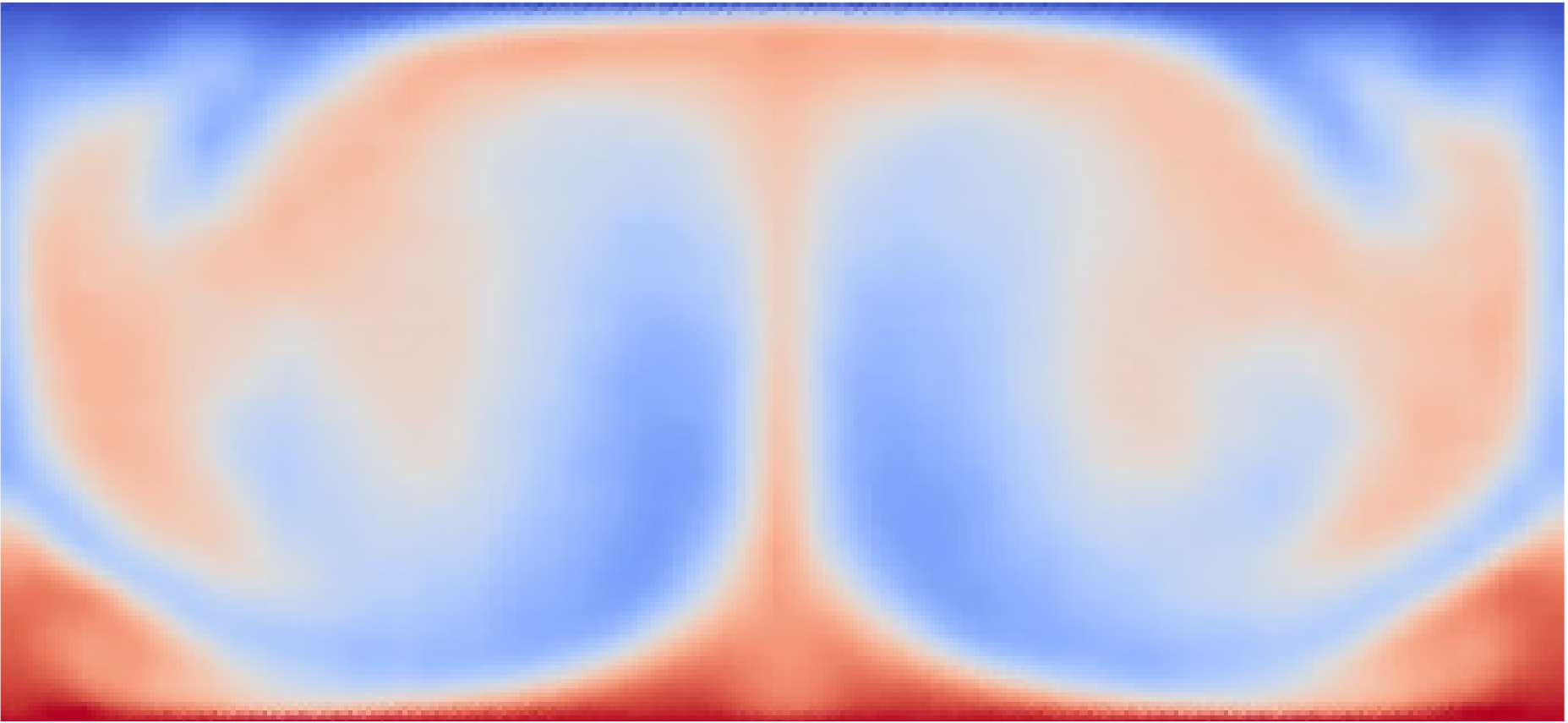}
		\caption{$\epsilon_{\tt CAE-PHODMD}=4.098\times 10^{-2}, n_{\tt delay}=16$}
	\end{subfigure}
	\caption{2D Rayleigh-B\'enard convection. $t=15.809$, $Ra=6\times 10^{6}$, with $n_{\tt latent}=10$.}
	\label{fig:2DRBC_u_t1_latent10}
\end{figure}

%

\begin{figure}[hbt!]
	\centering
	\begin{subfigure}[b]{0.45\textwidth}
		\centering
		\includegraphics[width=1.0\textwidth]{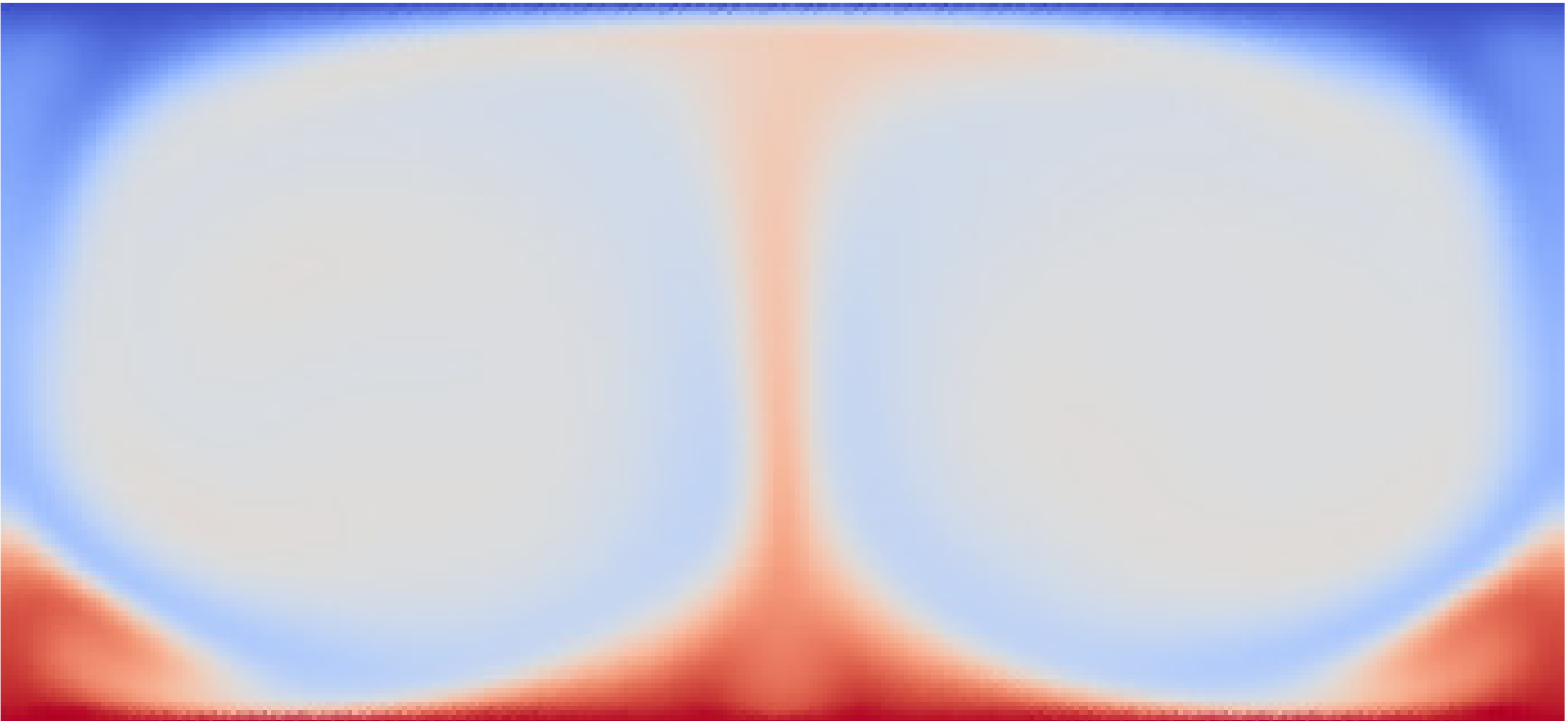}
		\caption{ground truth}
	\end{subfigure}
	\begin{subfigure}[b]{0.45\textwidth}
		\centering
		\includegraphics[width=1.0\textwidth]{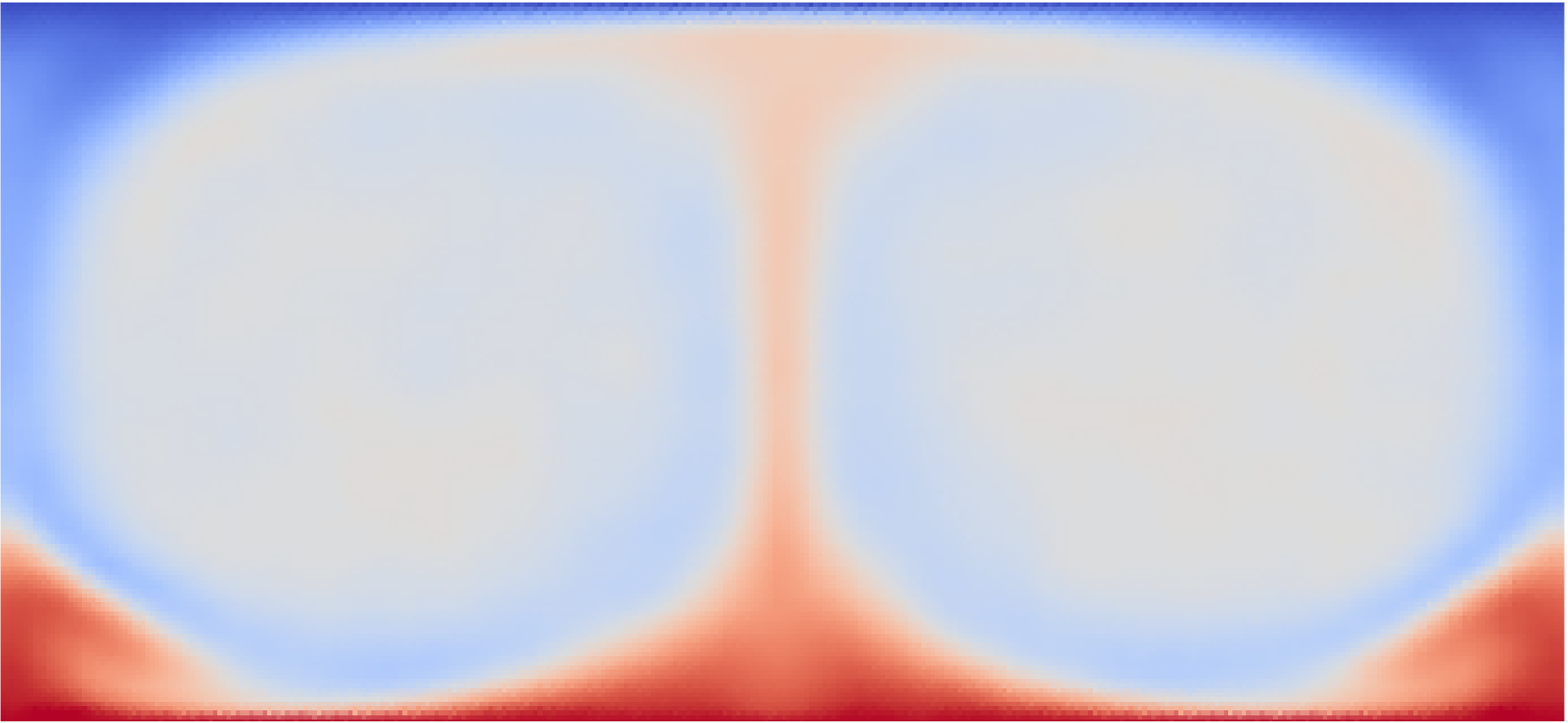}
		\caption{$\epsilon_{\tt CAE}=2.764\times 10^{-2}$}
	\end{subfigure}
	
	\begin{subfigure}[b]{0.45\textwidth}
		\centering
		\includegraphics[width=1.0\textwidth]{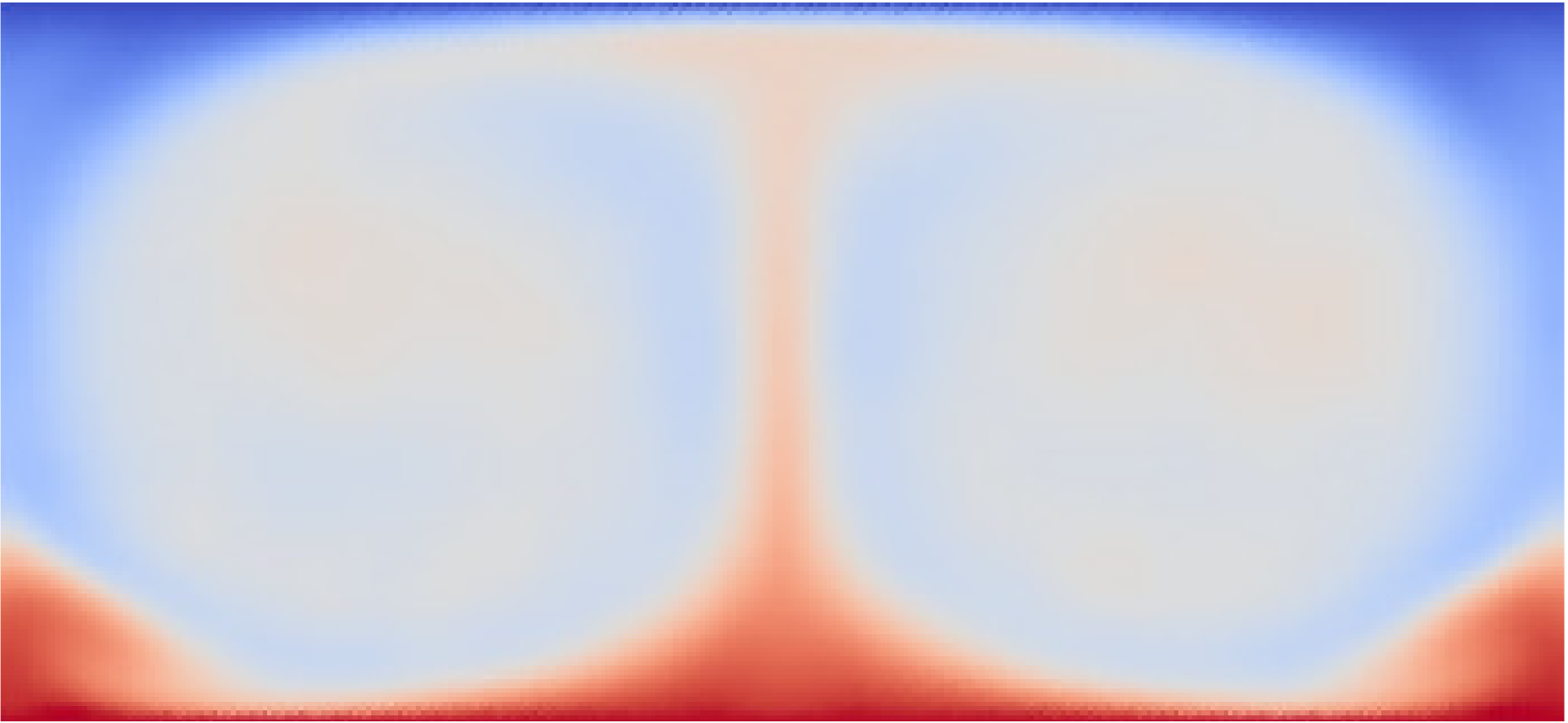}
		\caption{$\epsilon_{\tt CAE-PHODMD}=6.319\times 10^{-2}, n_{\tt delay}=4$}
	\end{subfigure}
	\begin{subfigure}[b]{0.45\textwidth}
		\centering
		\includegraphics[width=1.0\textwidth]{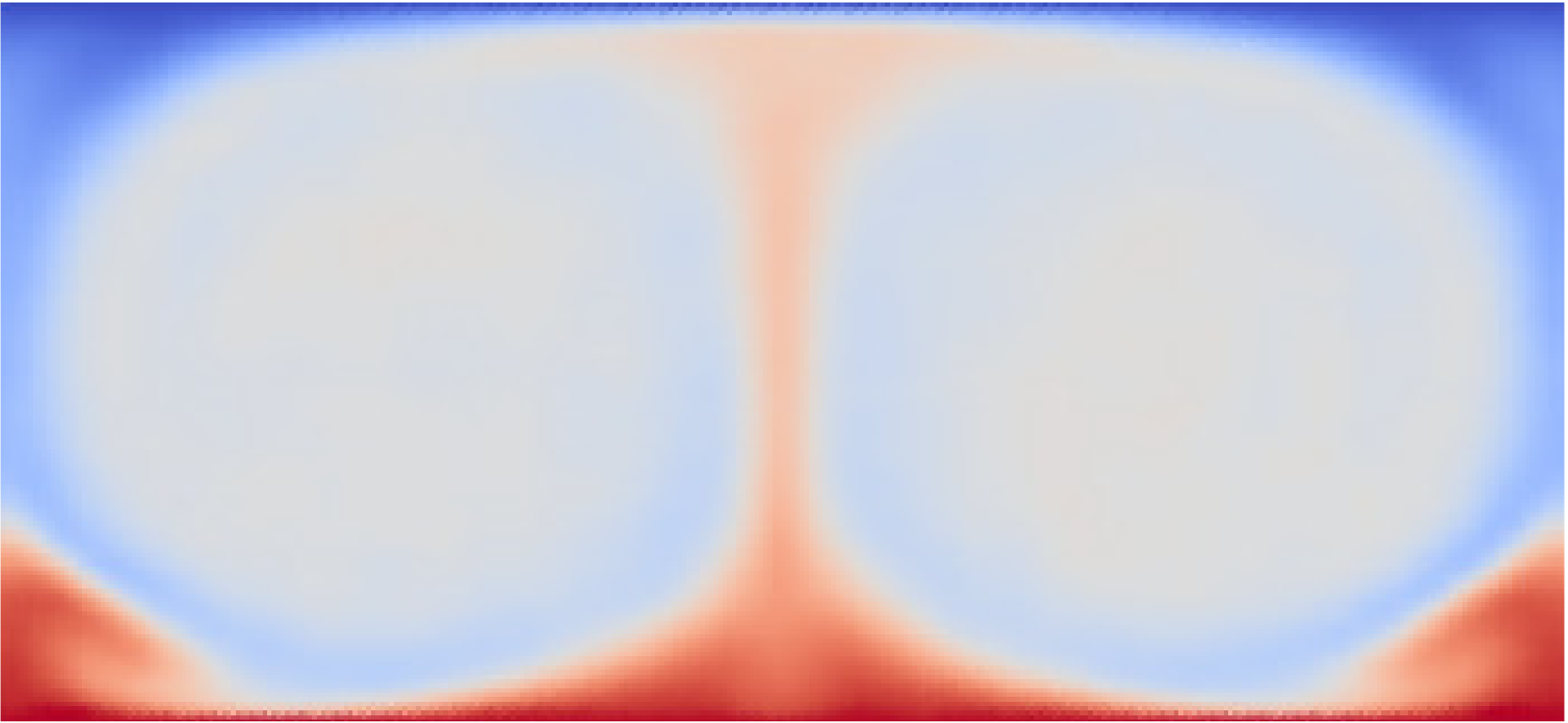}
		\caption{$\epsilon_{\tt CAE-PHODMD}=2.596\times 10^{-2}, n_{\tt delay}=16$}
	\end{subfigure}
	\caption{2D Rayleigh-B\'enard convection. $t=37.139$, $Ra=6\times 10^{6}$, with $n_{\tt latent}=10$.}
	\label{fig:2DRBC_u_t2_latent10}
\end{figure}

\begin{figure}[hbt!]
	\centering
	\begin{subfigure}[b]{0.45\textwidth}
		\centering
		\includegraphics[width=1.0\textwidth]{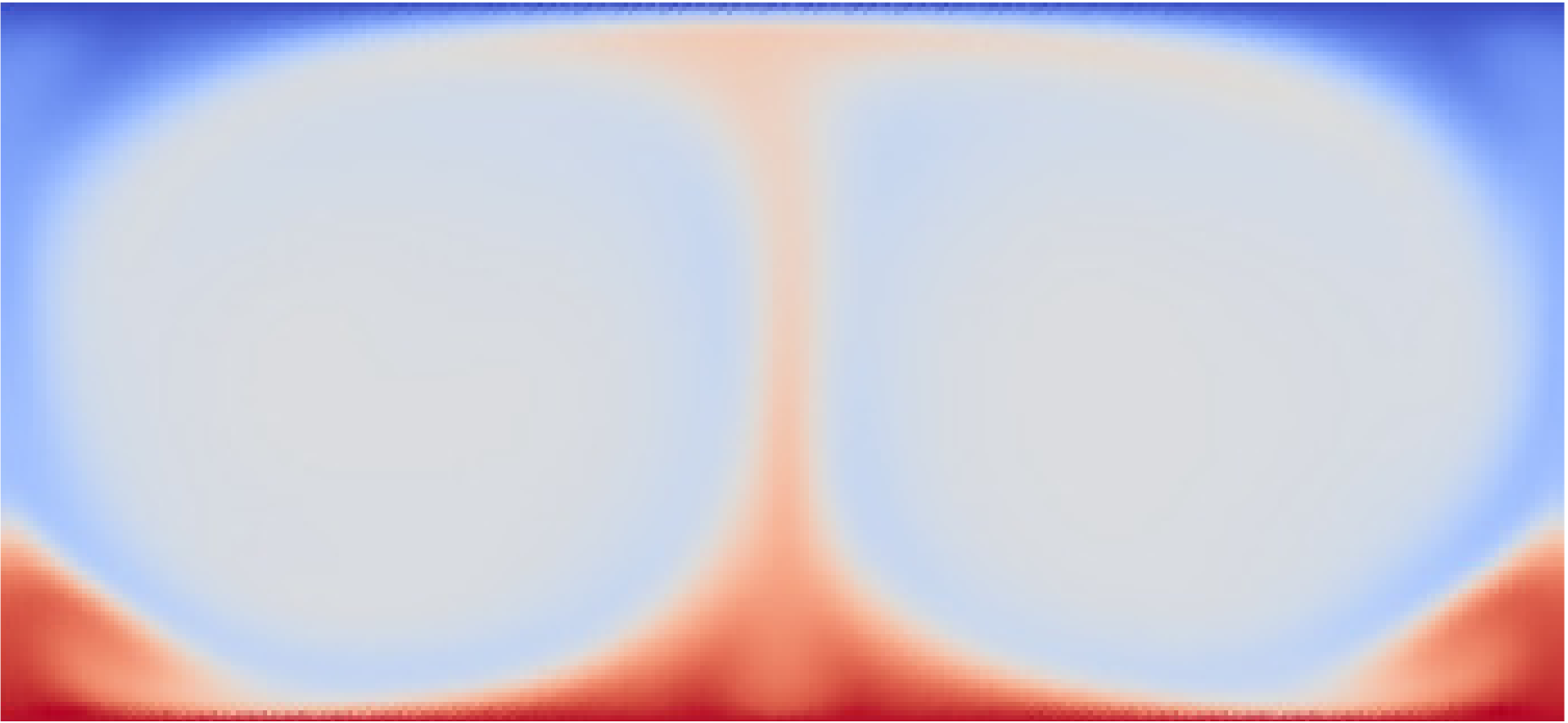}
		\caption{ground truth}
	\end{subfigure}
	\begin{subfigure}[b]{0.45\textwidth}
		\centering
		\includegraphics[width=1.0\textwidth]{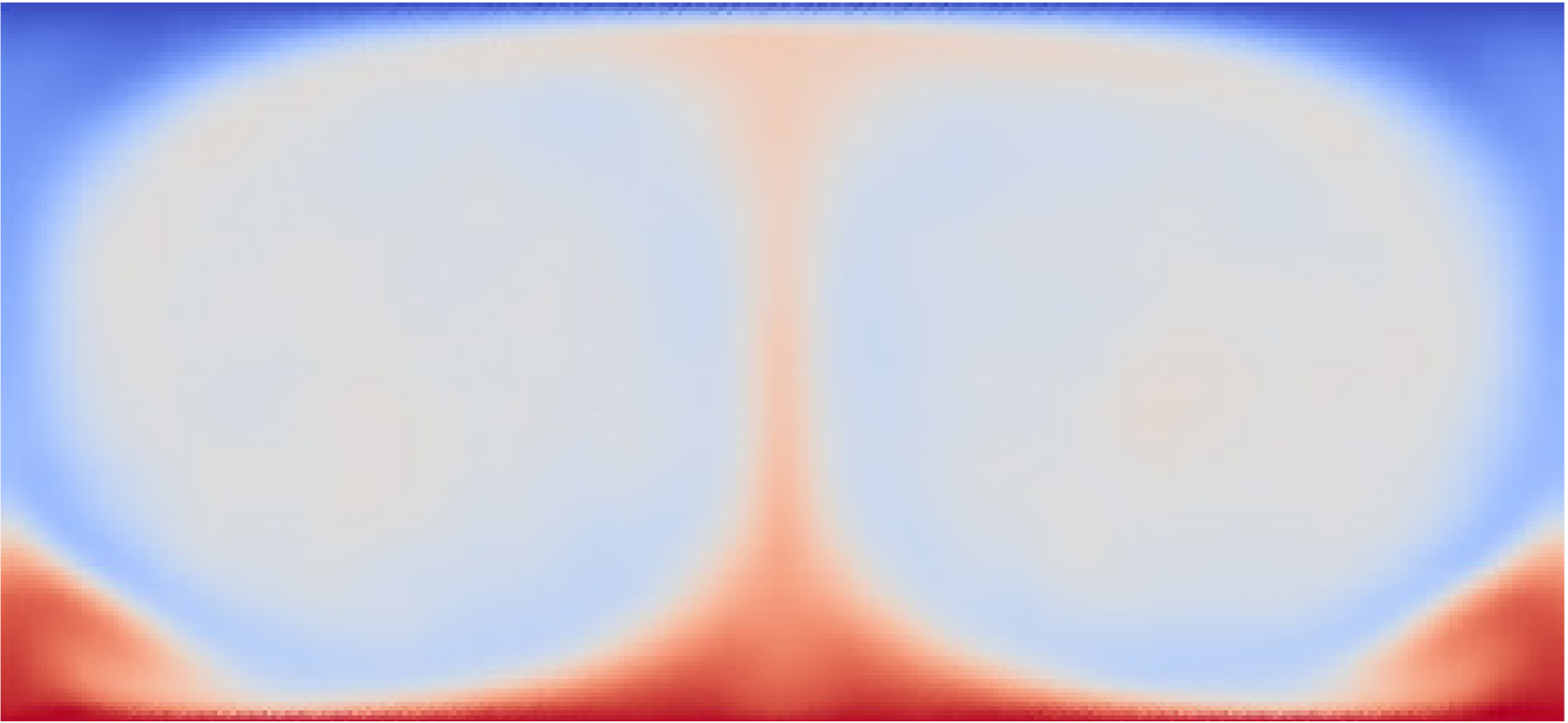}
		\caption{$\epsilon_{\tt CAE}=3.304\times 10^{-2}$}
	\end{subfigure}
	
	\begin{subfigure}[b]{0.45\textwidth}
		\centering
		\includegraphics[width=1.0\textwidth]{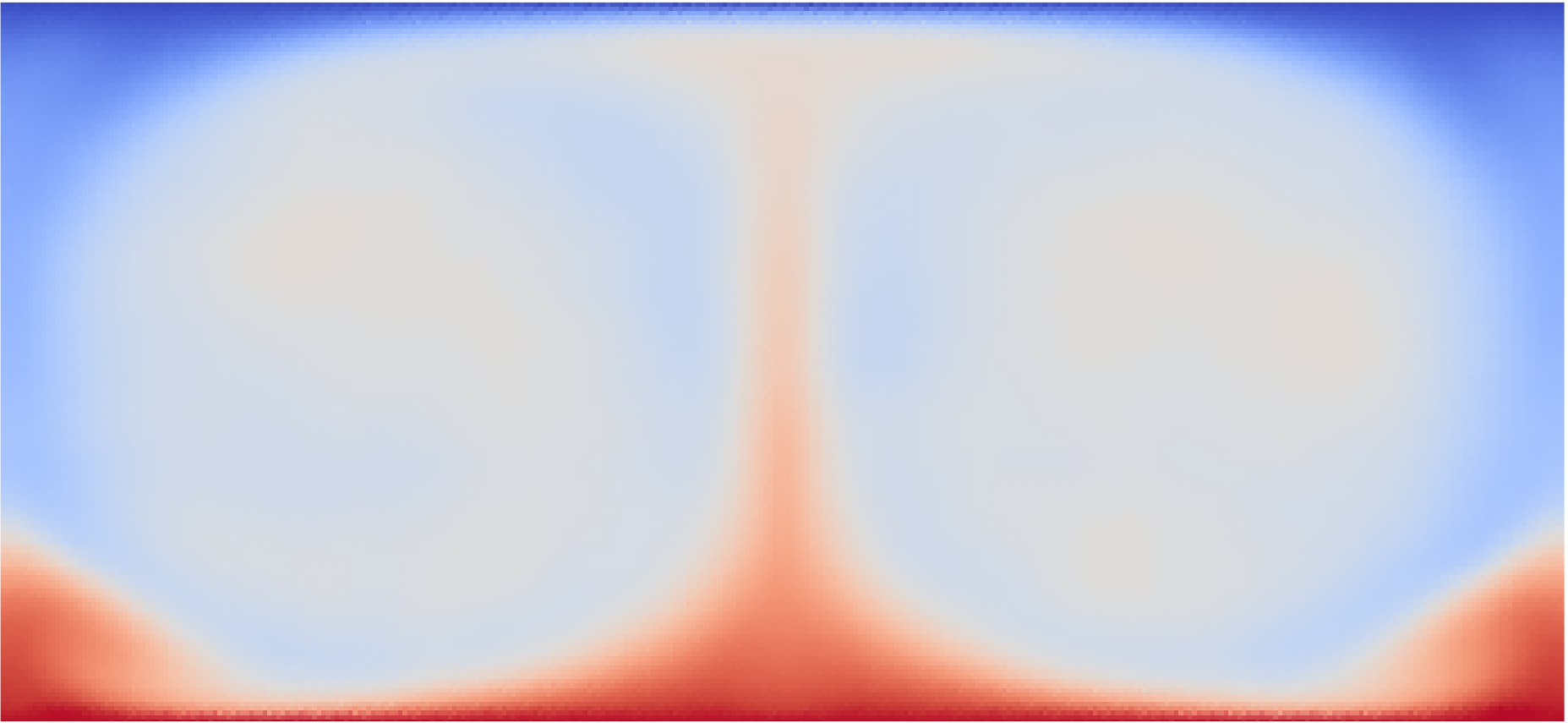}
		\caption{$\epsilon_{\tt CAE-PHODMD}=6.534\times 10^{-2}, n_{\tt delay}=4$}
	\end{subfigure}
	\begin{subfigure}[b]{0.45\textwidth}
		\centering
		\includegraphics[width=1.0\textwidth]{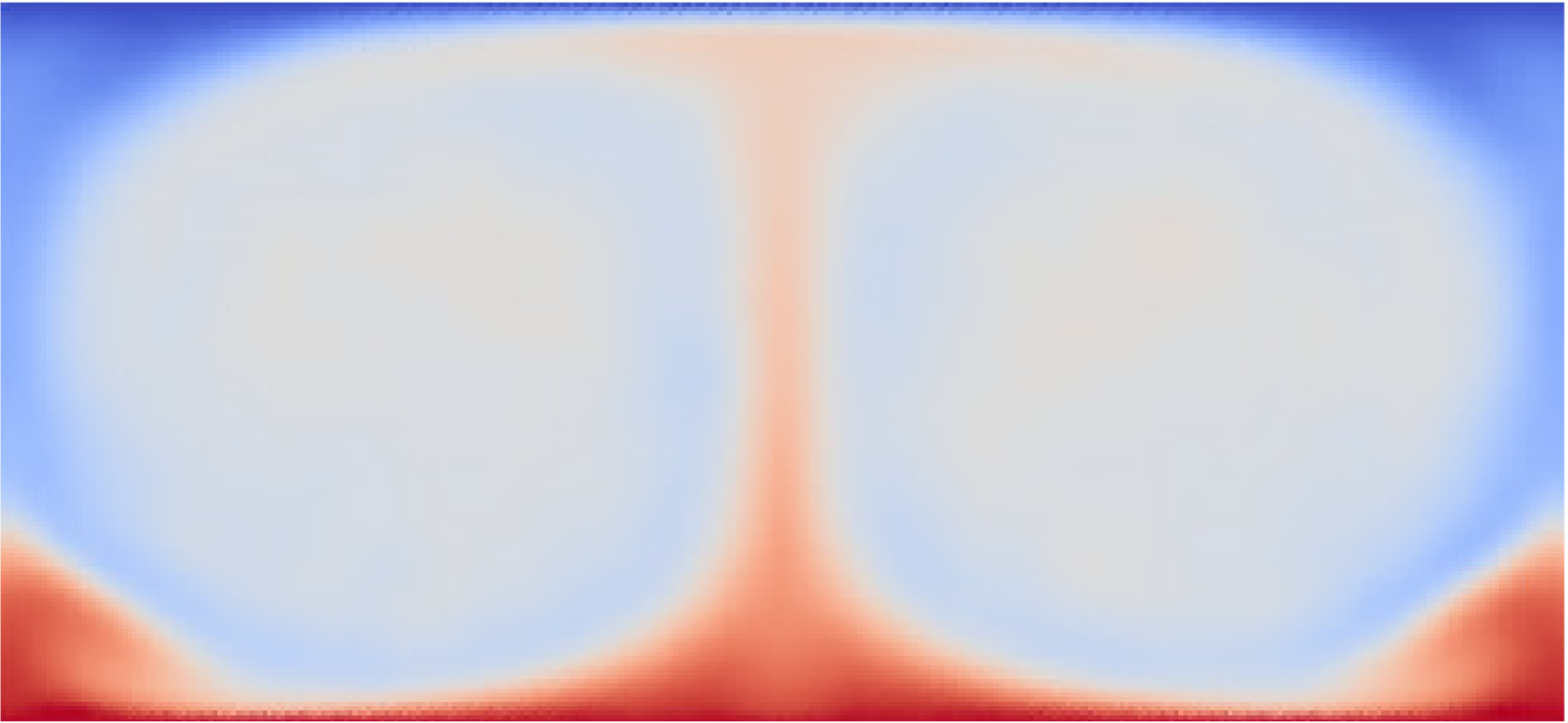}
		\caption{$\epsilon_{\tt CAE-PHODMD}=3.716\times 10^{-2}, n_{\tt delay}=8$}
	\end{subfigure}
	\caption{2D Rayleigh-B\'enard convection. $t=47.656$, $Ra=6\times 10^{6}$, with $n_{\tt latent}=10$.}
	\label{fig:2DRBC_u_t3_latent10}
\end{figure}

The ground truth and the corresponding solutions obtained by the CAE and {\tt CAE-PHODMD}
for the testing parameter value $Ra=6\times 10^{6}$ and $t=15.809$ with $n_{\tt latent}=10$
and $n_{\tt delay}=4, 16$ are shown in Figure \ref{fig:2DRBC_u_t1_latent10},
in which the errors are also listed.
The reconstruction by the CAE is accurate and captures small scale features.
For $n_{\tt latent}=10$, $n_{\tt delay}=4$ is not enough to ensure accurate prediction,
as the important patterns are lost.
The result with $n_{\tt delay}=16$ is nearly as accurate as the one obtained by the CAE.
The solutions at $t=37.139$ are shown in Figure \ref{fig:2DRBC_u_t2_latent10}.
It is observed that with $n_{\tt latent}=10$, increasing the number of time-delay embedding reduces the error,
and $n_{\tt delay}=16$ suffices to recover accurate results.
In this test, we also check the extrapolation capability by the CAE and {\tt CAE-PHODMD}
at $t=47.656$ with $n_{\tt latent}=10$, given in Figures \ref{fig:2DRBC_u_t3_latent10}.
The reconstruction by the CAE is stable and the result is accurate.
The result with $n_{\tt latent}=10$ and $n_{\tt delay}=8$ captures the main features accurately,
and is comparable to that obtained by the CAE,
showing the prediction ability of the approach.


%

\begin{figure}[hbt!]
	\centering
	\includegraphics[width=1.0\textwidth]{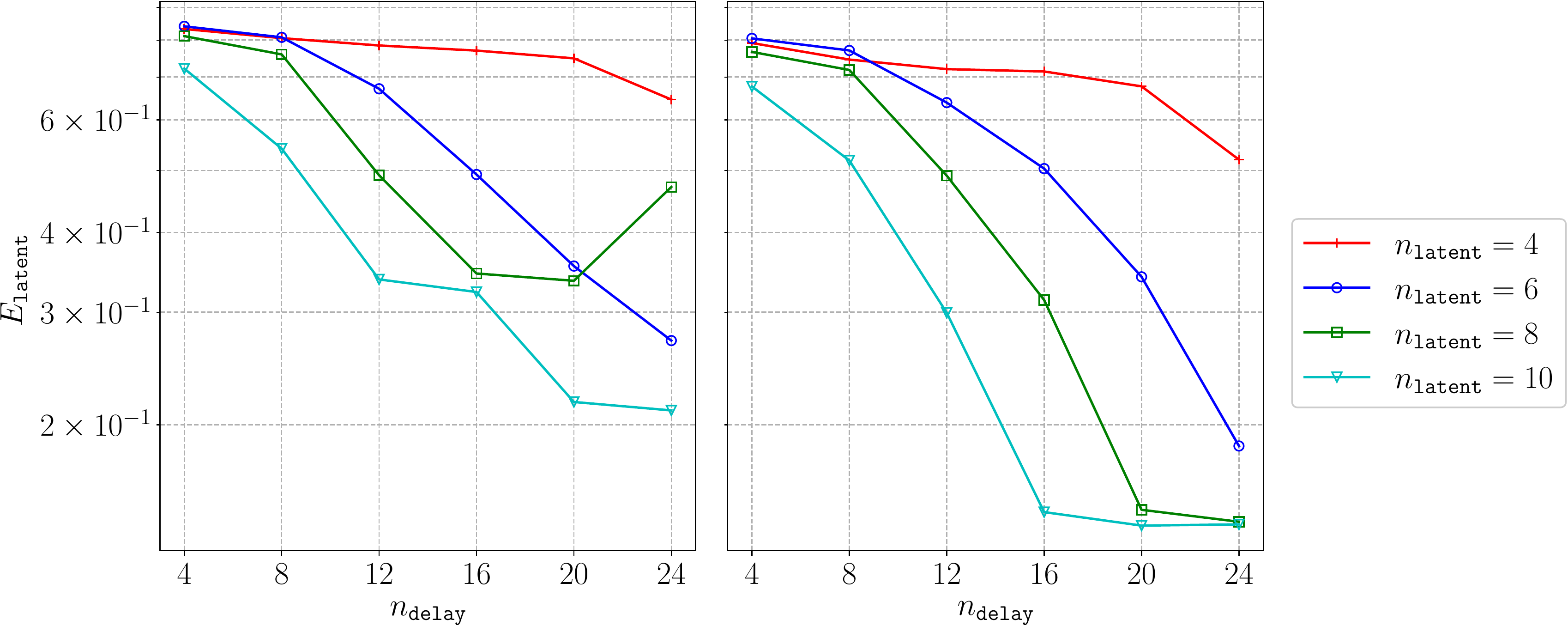}
	\caption{2D Rayleigh-B\'enard convection.
	The errors in the latent space $E_{\tt latent}$ w.r.t. the number of the time-delay embedding $n_{\tt delay}$ for different dimensions of the latent space $n_{\tt latent}$.
	The left figure is for all testing times, while the right for interpolation in time.}
	\label{fig:2DRBC_err_q_delay}
\end{figure}

\begin{figure}[hbt!]
	\centering
	\includegraphics[width=0.6\textwidth]{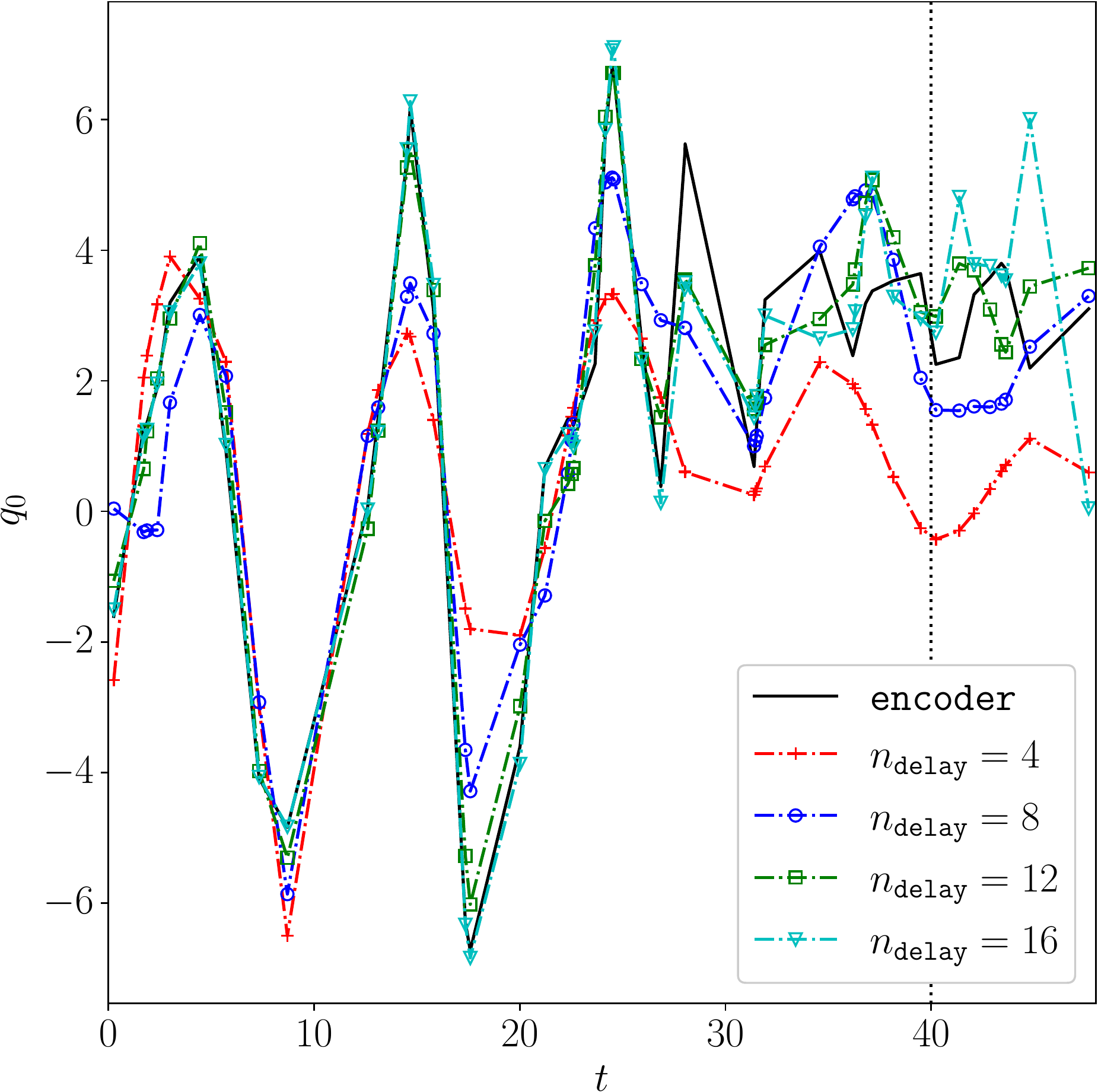}
	\caption{2D Rayleigh-B\'enard convection.
		The temporal evolution of the first latent variable for $Ra=6\times 10^{6}$ with $n_{\tt latent}=10$ and different $n_{\tt delay}$.}
	\label{fig:2DRBC_q}
\end{figure}

\begin{figure}[hbt!]
	\centering
	\includegraphics[width=0.6\textwidth]{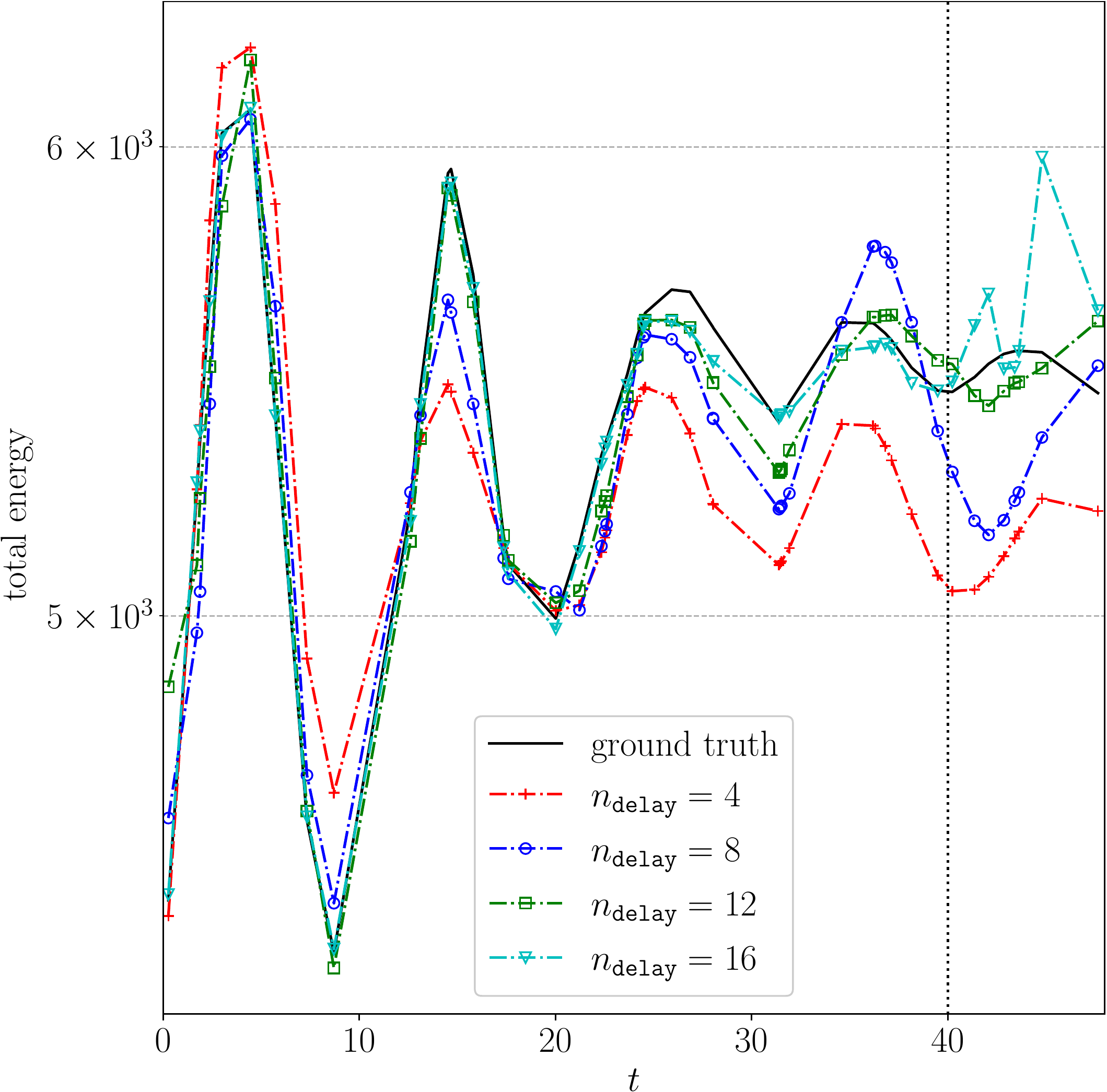}
	\caption{2D Rayleigh-B\'enard convection.
		The temporal evolution of the total energy for $Ra=6\times 10^{6}$ with $n_{\tt latent}=10$ and different $n_{\tt delay}$.}
	\label{fig:2DRBC_energy}
\end{figure}

Figure \ref{fig:2DRBC_err_q_delay} plots the errors $E_{\tt latent}$ with respect to $n_{\tt delay}$ for different $n_{\tt latent}$.
The left and right panels correspond to all testing times and interpolation in time, respectively.
It is expected that the errors in the left including extrapolation errors are larger,
and that the interpolation errors decrease as $n_{\tt delay}$ increases.
To further understand the dynamics in the latent space, Figure \ref{fig:2DRBC_q} shows the evolution of the first latent variable for $n_{\tt latent}=10$ with different $n_{\tt delay}$.
One can see that before $t=40$, i.e., for interpolation, the parametric HODMD gets more accurate as $n_{\tt delay}$ increases,
and the errors are slightly larger after $t=40$.
This can also be observed from the evolution of the total energy in Figure \ref{fig:2DRBC_energy}.


\subsection{2D Kelvin-Helmholtz instability}\label{sec:2DKHI}
\subsubsection{Setup}
In this section, the 2D Kelvin-Helmholtz instability which is a transport-dominated problem is considered.
The full-order solutions are obtained by solving the following compressible Euler equations
\begin{equation*}
	\dfrac{\partial}{\partial t}
	\begin{pmatrix}
		\rho \\ \rho v_x \\ \rho v_y \\ E \\
	\end{pmatrix}
	+ \dfrac{\partial}{\partial x}
	\begin{pmatrix}
		\rho v_x \\ \rho v_x^2+p \\  \rho v_x v_y \\  (E+p)v_x \\
	\end{pmatrix}
	+ \dfrac{\partial}{\partial y}
	\begin{pmatrix}
		\rho v_y \\ \rho v_x v_y \\  \rho v_y^2+p \\  (E+p)v_y \\
	\end{pmatrix} = 0,
\end{equation*}
where $\rho, v_x, v_y, E$ are the mass density, velocities, and total energy, respectively.
The pressure $p$ can be obtained by the perfect gas equation of state
\begin{equation*}
	p = (E - \frac12\rho(v_x^2+v_y^2))/(\Gamma-1),
\end{equation*}
with the adiabatic index $\Gamma=1.4$ in this test.
The initial data are
\begin{align*}
	&y_1 = 0.25, ~y_2 = 0.75,\\
	&\begin{cases}
		\rho=2,~ v_x=0.5+\omega, &\text{if}~ y_1 < y < y_2, \\
		\rho=1,~ v_x=0.5-\omega, &\text{else}, \\
	\end{cases}\\
	&v_y = 0.1\sin(4\pi x)(\exp(-0.5(y-y_1)^2)/0.05 + \exp(-0.5(y-y_2)^2)/0.05), \\
	&p = 2.5,
\end{align*}
where $\omega$ controls the shear velocity of the initial discontinuities.
A fifth-order finite difference WENO scheme \cite{Jiang1996Efficient} is employed to obtain the full-order solutions on the uniform $128\times 128$ mesh.
The datasets are chosen as follows.
The training set comprises uniform sampling in the parameter space ranging from $-0.1$ to $0.1$ with step size $0.01$,
except for $\{-0.03,0.08,-0.04,0.05\}$,
where the first two is for validation and the last two for testing.
The training set comprises $100$ times uniform in $[0.01,1]$ with $\Delta t=0.01$,
the validation and testing sets comprise $20$ and $10$ random times in $[0.01,1]$, respectively.
The kernel size in all the convolutional and deconvolutional layers is $5$ with the stride and padding as $2$.

\subsubsection{Results}
We perform a grid search with different weight decay $10^{-8}$, $10^{-9}$, $10^{-10}$, $10^{-11}$,
number of convolutional layers $4$, $5$, $6$, and the dimension of the latent space $n_{\tt latent}=6$, $8$, $10$, $15$.
The final CAE architecture is shown in Table \ref{tab:2DKHI_arch}.

\begin{table}[hbt!]
	\centering
	\begin{tabular}{c|r|r}
		\toprule
		\multicolumn{3}{c}{encoder} \\ \hline
		layer & input shape & output shape \\ \hline
		{\tt Conv2d with SiLU} & $4\times 128\times 128$ & $16\times 64\times 64$ \\
		{\tt Conv2d with SiLU} & $16\times 64\times 64$ & $16\times 32\times 32$ \\
		{\tt Conv2d with SiLU} & $16\times 32\times 32$ & $16\times 16\times 16$ \\
		{\tt Conv2d with SiLU} & $16\times 16\times 16$ & $16\times 8\times 8$ \\
		{\tt Flatten} & $16\times 8\times 8$ & $ 1024 $ \\
		{\tt Linear with SiLU} & $ 1024 $ & $90$ \\
		{\tt Linear with SiLU} & $ 90 $ & $8$ \\
		\toprule
		\multicolumn{3}{c}{decoder} \\ \hline
		layer & input shape & output shape \\ \hline		
		{\tt Linear with SiLU} & $8$ & $90$ \\
		{\tt Linear with SiLU} & $90$ & $1024$ \\
		{\tt Unflatten} & $ 1024 $ & $16\times 8\times 8$ \\
		{\tt ConvTranspose2d with SiLU} & $16\times 8\times 8$ & $16\times 16\times 16$ \\
		{\tt ConvTranspose2d with SiLU} & $16\times 16\times 16$ & $16\times 32\times 32$ \\
		{\tt ConvTranspose2d with SiLU} & $16\times 32\times 32$ & $16\times 64\times 64$ \\
		{\tt ConvTranspose2d with SiLU} & $16\times 64\times 64$ & $4\times 128\times 128$ \\
		\bottomrule
	\end{tabular}
	\caption{2D Kelvin-Helmholtz instability: The architecture of the CAE during the grid search with the best reconstruction error on the validation set.
	The weight decay is $10^{-8}$ in the training.}
	\label{tab:2DKHI_arch}
\end{table}

%

\begin{figure}[hbt!]
	\centering
	\includegraphics[width=1.0\textwidth]{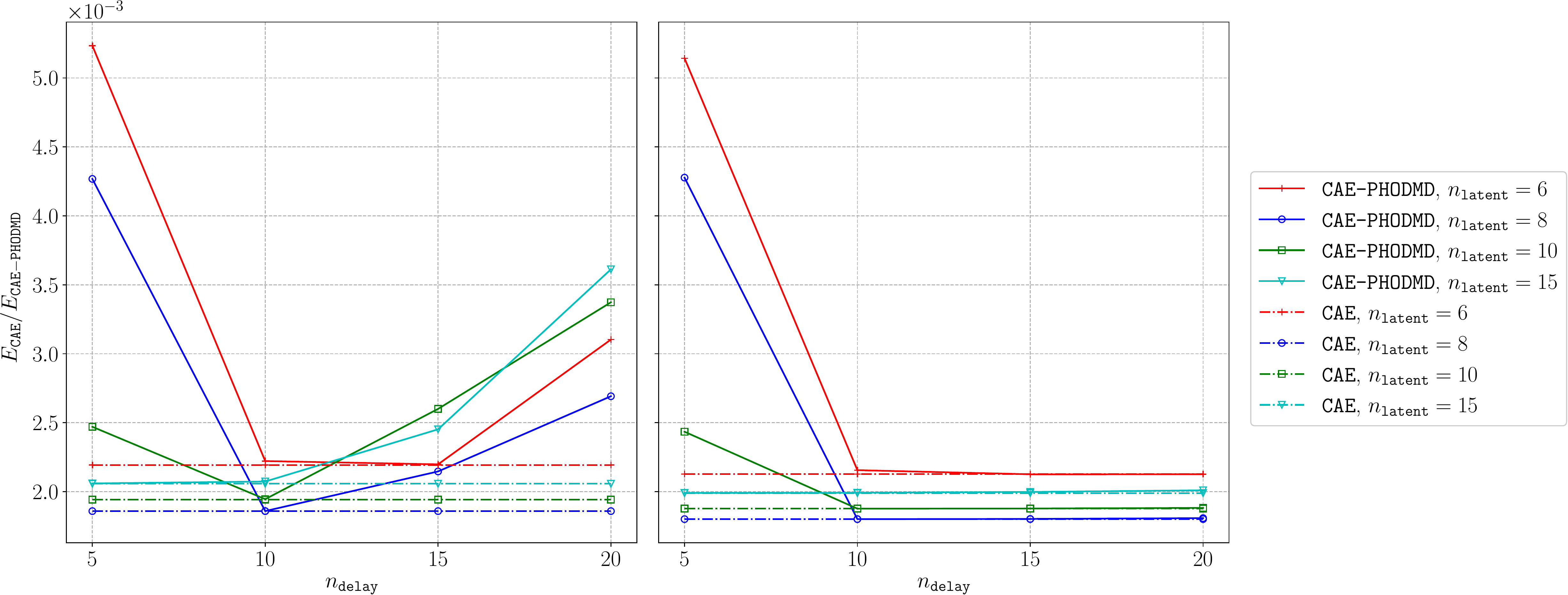}
	\caption{2D Kelvin-Helmholtz instability. The reconstruction errors $E_{\tt CAE}$ and prediction errors $E_{\tt CAE-PHODMD}$ w.r.t. the number of the time-delay embedding $n_{\tt delay}$ for different dimensions of the latent space $n_{\tt latent}$.
	The left figure is for all the testing times, while the right excludes one trailing time.}
	\label{fig:2DKHI_err_u_CAE_PHODMD}
\end{figure}

The errors $E_{\tt CAE}$ and $E_{\tt CAE-PHODMD}$ with different $n_{\tt latent}$ and $n_{\tt delay}$ are given in Figure \ref{fig:2DKHI_err_u_CAE_PHODMD}.
The left figure considers all the samplings in the testing set, while the right excludes one trailing time outside $[0,0.9]$.
It can be seen that the results of {\tt CAE-PHODMD} are accurate away from the last testing time.
For a given $n_{\tt latent}$, the errors $E_{\tt CAE-PHODMD}$ in the right figure decrease as $n_{\tt delay}$ increases,
then remain constant and are comparable to the reconstruction errors obtained by the CAE.
One can conclude that using $n_{\tt latent}=8$ with $n_{\tt delay}=10$ is able to recover accurate results.

\begin{figure}[hbt!]
\centering
\begin{subfigure}[b]{0.45\textwidth}
	\centering
	\includegraphics[width=1.0\textwidth]{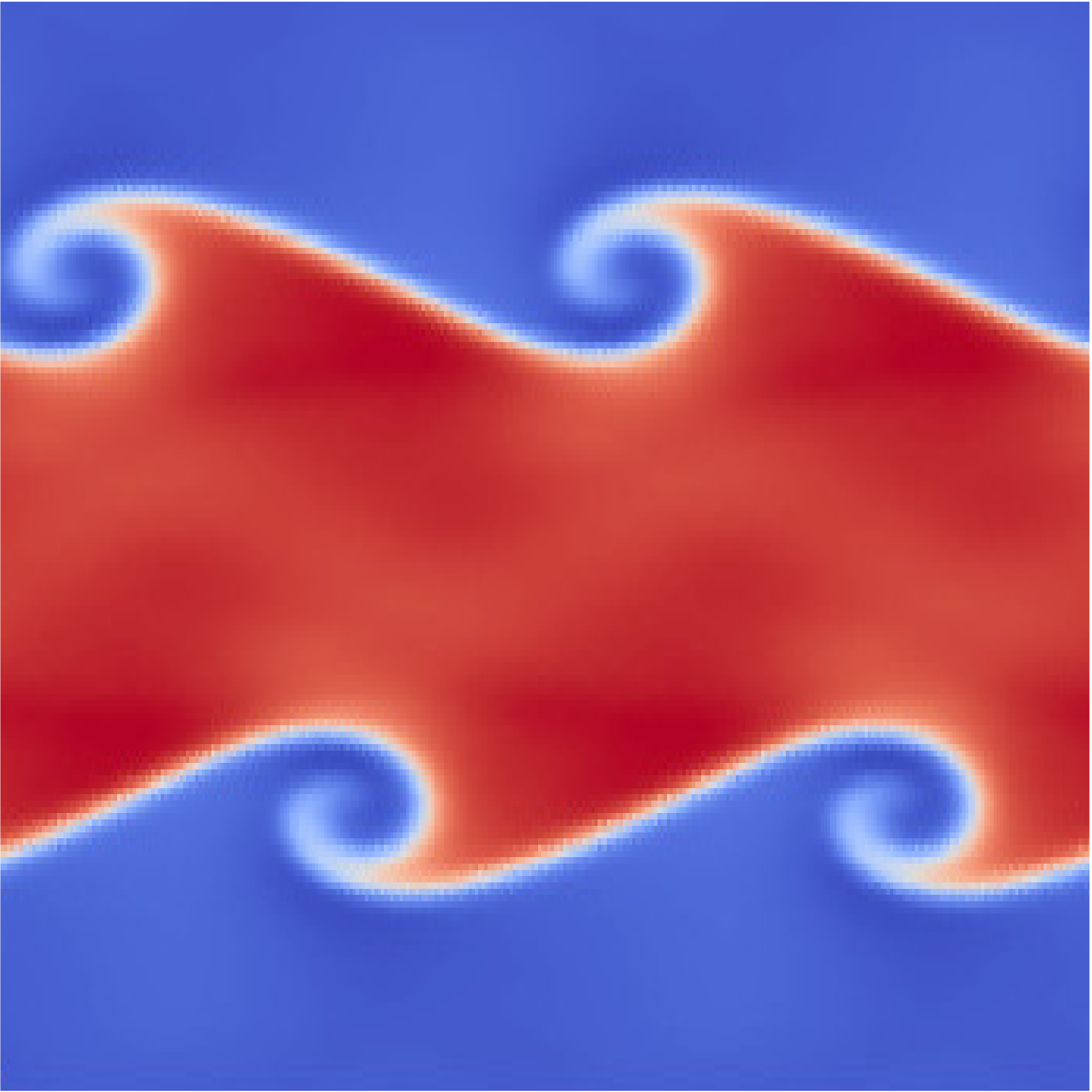}
	\caption{ground truth}
\end{subfigure}
\begin{subfigure}[b]{0.45\textwidth}
	\centering
	\includegraphics[width=1.0\textwidth]{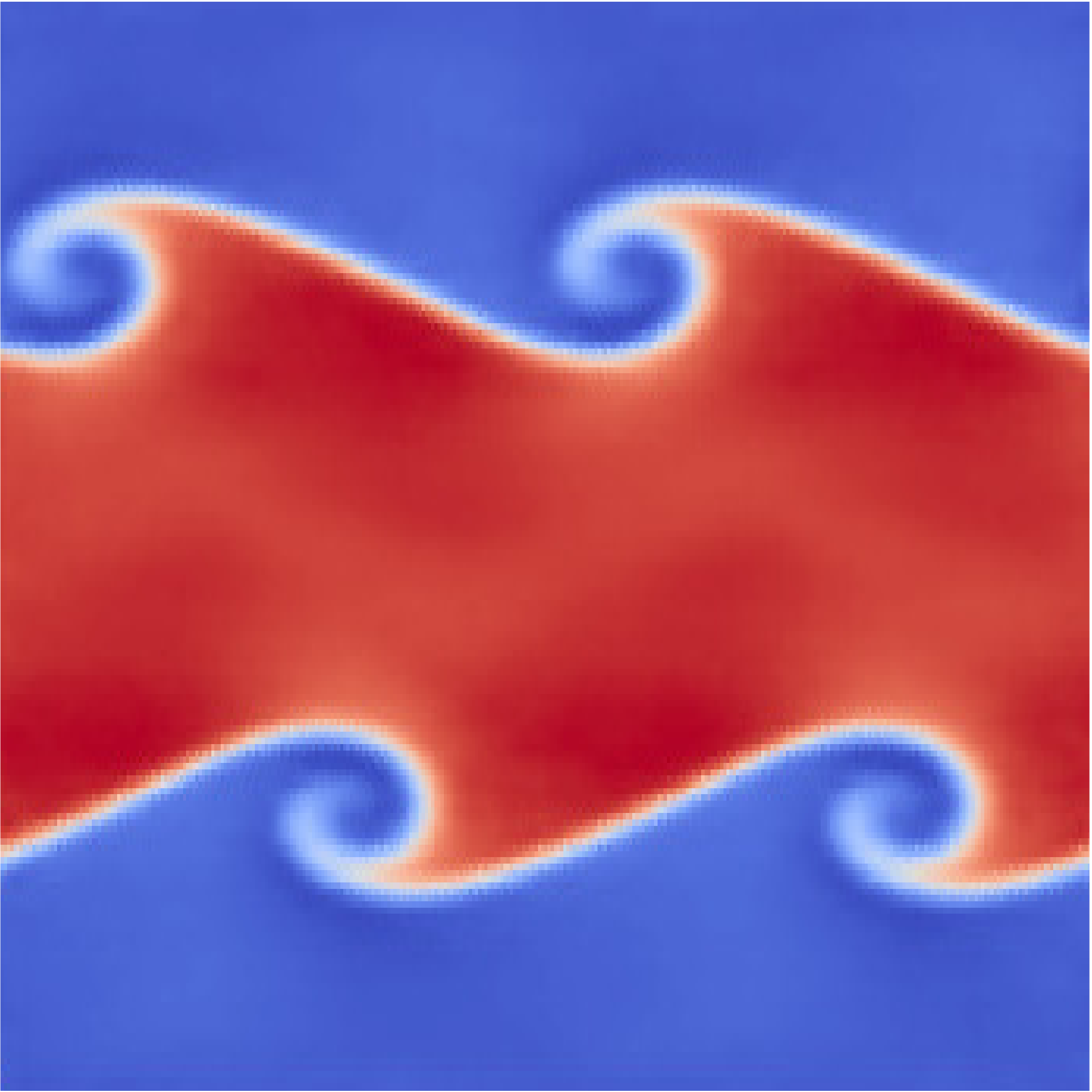}
	\caption{$\epsilon_{\tt CAE}=4.090\times 10^{-3}$}
\end{subfigure}

\begin{subfigure}[b]{0.45\textwidth}
	\centering
	\includegraphics[width=1.0\textwidth]{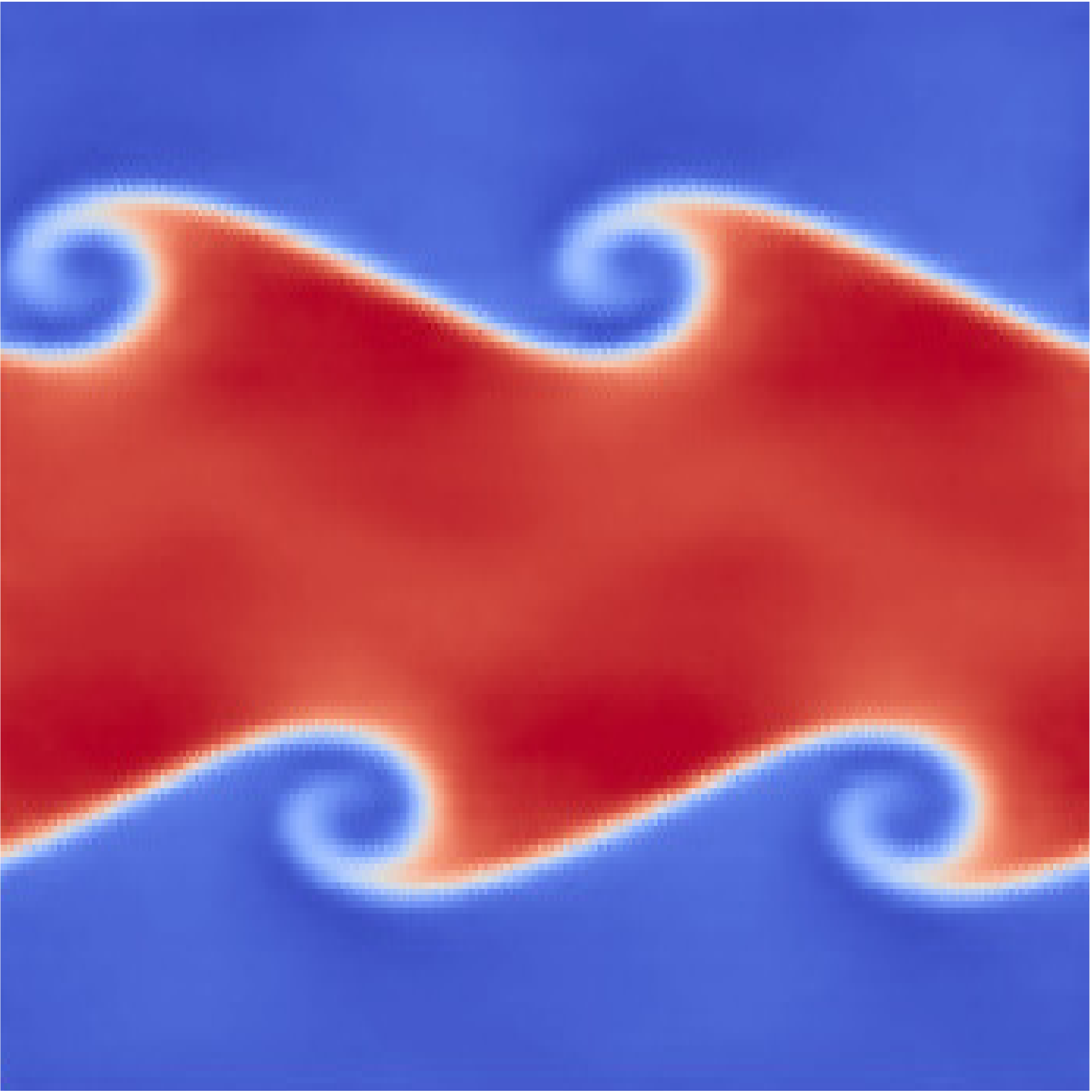}
	\caption{$\epsilon_{\tt CAE-PHODMD}=9.921\times 10^{-3}$, $n_{\tt delay}=5$}
\end{subfigure}
\begin{subfigure}[b]{0.45\textwidth}
	\centering
	\includegraphics[width=1.0\textwidth]{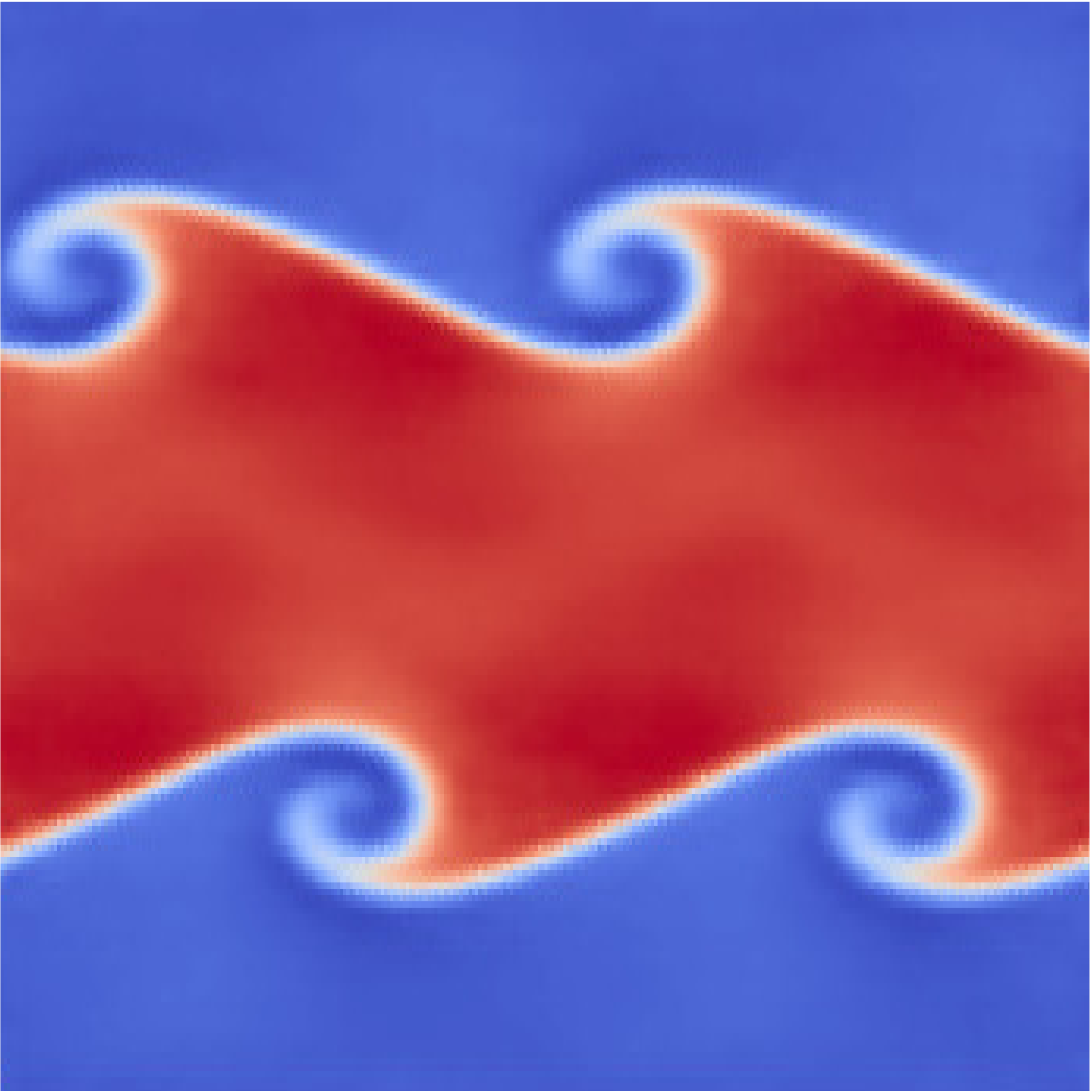}
	\caption{$\epsilon_{\tt CAE-PHODMD}=4.076\times 10^{-3}$, $n_{\tt delay}=10$}
\end{subfigure}
\caption{2D Kelvin-Helmholtz instability. $t=0.568$, $\omega=0.05$, with $n_{\tt latent}=8$.}
\label{fig:2DKHI_u_t1_w2_latent8}
\end{figure}

\begin{figure}[hbt!]
	\centering
	\begin{subfigure}[b]{0.45\textwidth}
		\centering
		\includegraphics[width=1.0\textwidth]{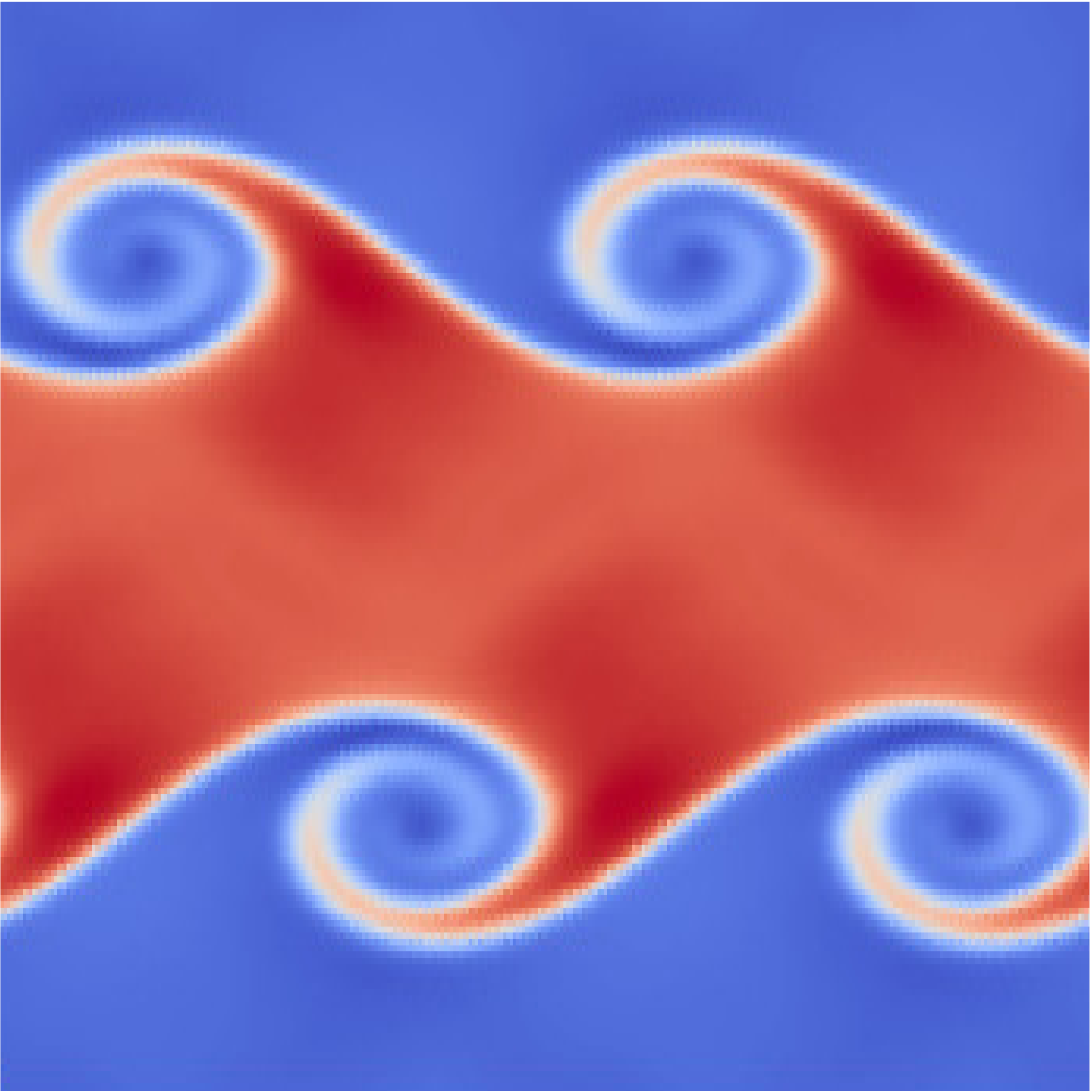}
		\caption{ground truth}
	\end{subfigure}
	\begin{subfigure}[b]{0.45\textwidth}
		\centering
		\includegraphics[width=1.0\textwidth]{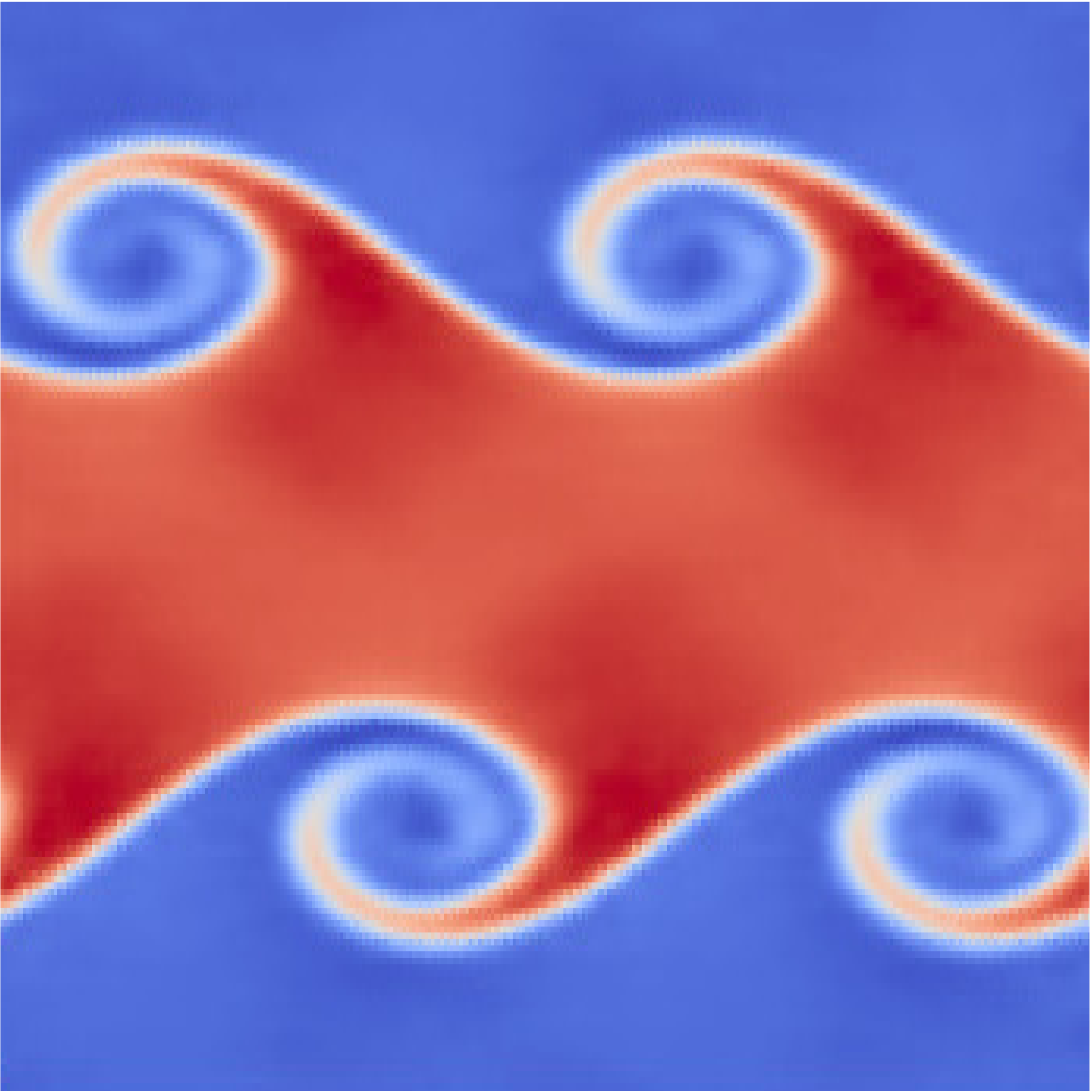}
		\caption{$\epsilon_{\tt CAE}=5.739\times 10^{-3}$}
	\end{subfigure}
	
	\begin{subfigure}[b]{0.45\textwidth}
		\centering
		\includegraphics[width=1.0\textwidth]{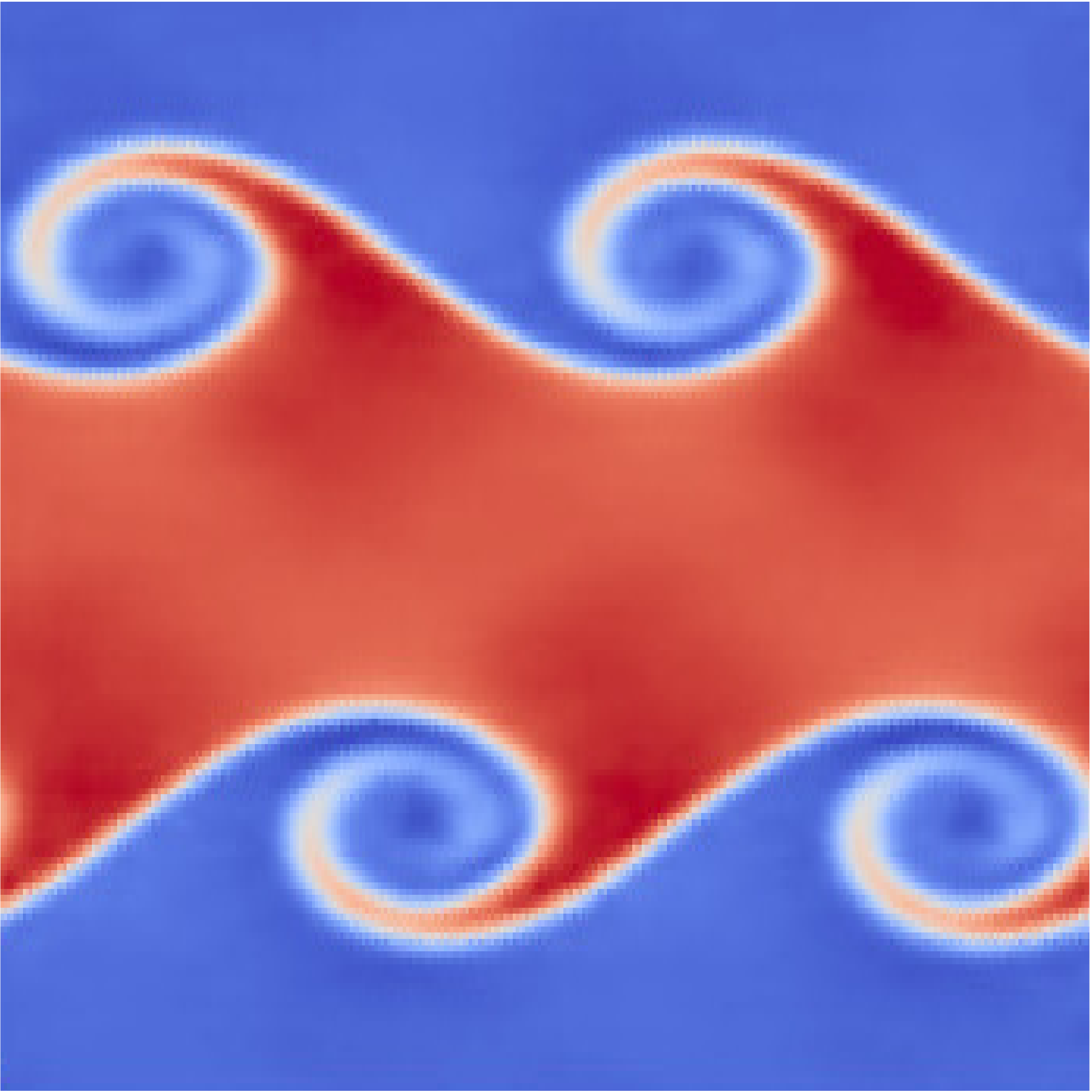}
		\caption{$\epsilon_{\tt CAE-PHODMD}=8.468\times 10^{-3}$, $n_{\tt delay}=5$}
	\end{subfigure}
	\begin{subfigure}[b]{0.45\textwidth}
		\centering
		\includegraphics[width=1.0\textwidth]{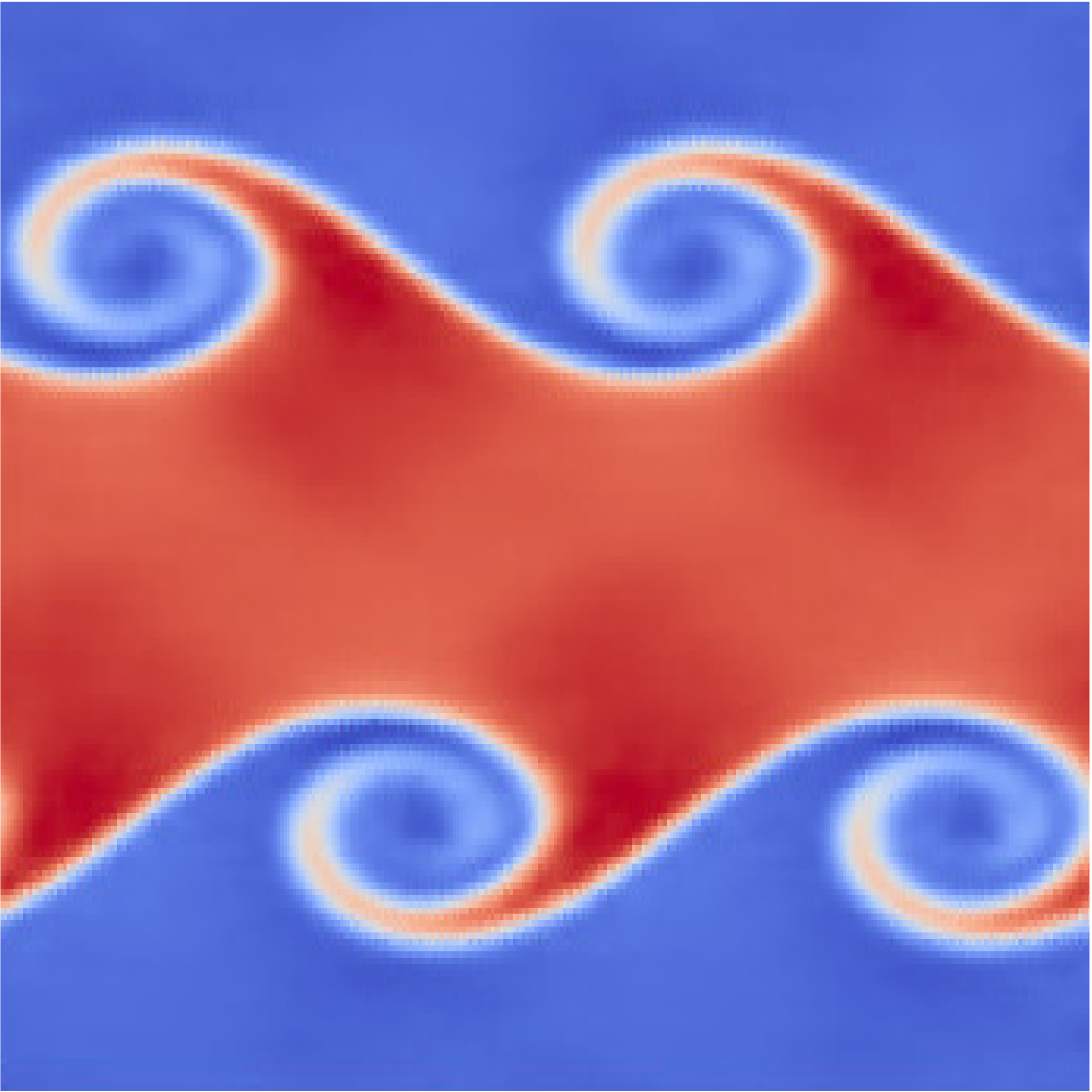}
		\caption{$\epsilon_{\tt CAE-PHODMD}=5.747\times 10^{-3}$, $n_{\tt delay}=10$}
	\end{subfigure}
	\caption{2D Kelvin-Helmholtz instability. $t=0.933$, $\omega=0.05$, with $n_{\tt latent}=8$.}
	\label{fig:2DKHI_u_t2_w2_latent8}
\end{figure}

The ground truth and corresponding solutions obtained by the CAE and {\tt CAE-PHODMD}
for a specific testing parameter value $w=0.05$ and time $t=0.568$ with $n_{\tt latent}=8$
and $n_{\tt delay}=5, 10$ are shown in Figure \ref{fig:2DKHI_u_t1_w2_latent8} with corresponding errors.
The reconstruction by the CAE is accurate and captures the rolls at the correct position.
For $n_{\tt latent}=8$ in this case,
using $n_{\tt delay}=5$ and $n_{\tt delay}=10$ both captures main features,
and the latter is more accurate,
which is expected since the HODMD uses more information based on a larger sliding window.
The results at $t=0.933$ are plotted in Figures \ref{fig:2DKHI_u_t2_w2_latent8}.
At this time, the patterns are more complex, and {\tt CAE-PHODMD} recovers accurate result
with $n_{\tt latent}=8$ and $n_{\tt delay}=10$, comparable to the reconstruction obtained by the CAE,
which shows that the approach captures the dynamics of this transport-dominated problem.

\begin{figure}[hbt!]
	\centering
	\includegraphics[width=1.0\textwidth]{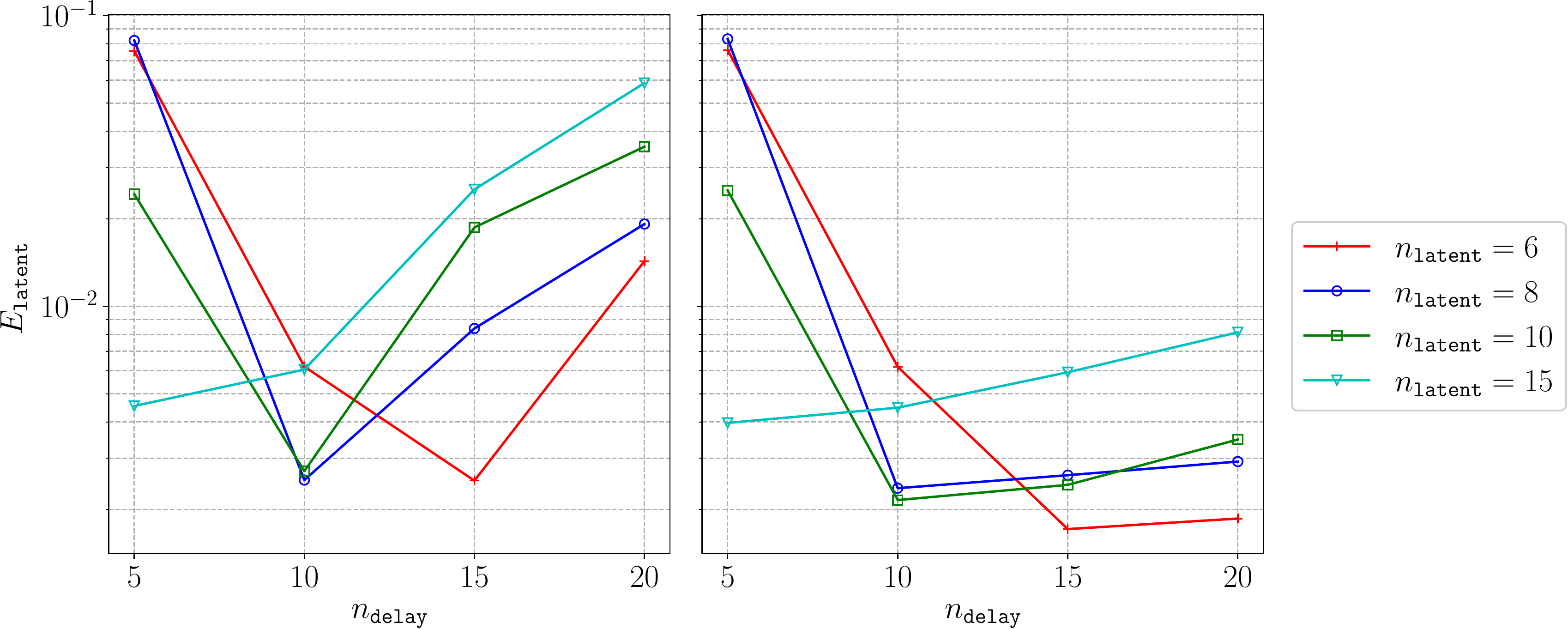}
	\caption{2D Kelvin-Helmholtz instability. The errors in the latent space $E_{\tt latent}$ w.r.t. the number of the time-delay embedding $n_{\tt delay}$ for different dimensions of the latent space $n_{\tt latent}$.
	The left figure is for all the testing times, while the right excludes one trailing time.}
	\label{fig:2DKHI_err_q_delay}
\end{figure}

\begin{figure}[hbt!]
	\centering
	\includegraphics[width=0.6\textwidth]{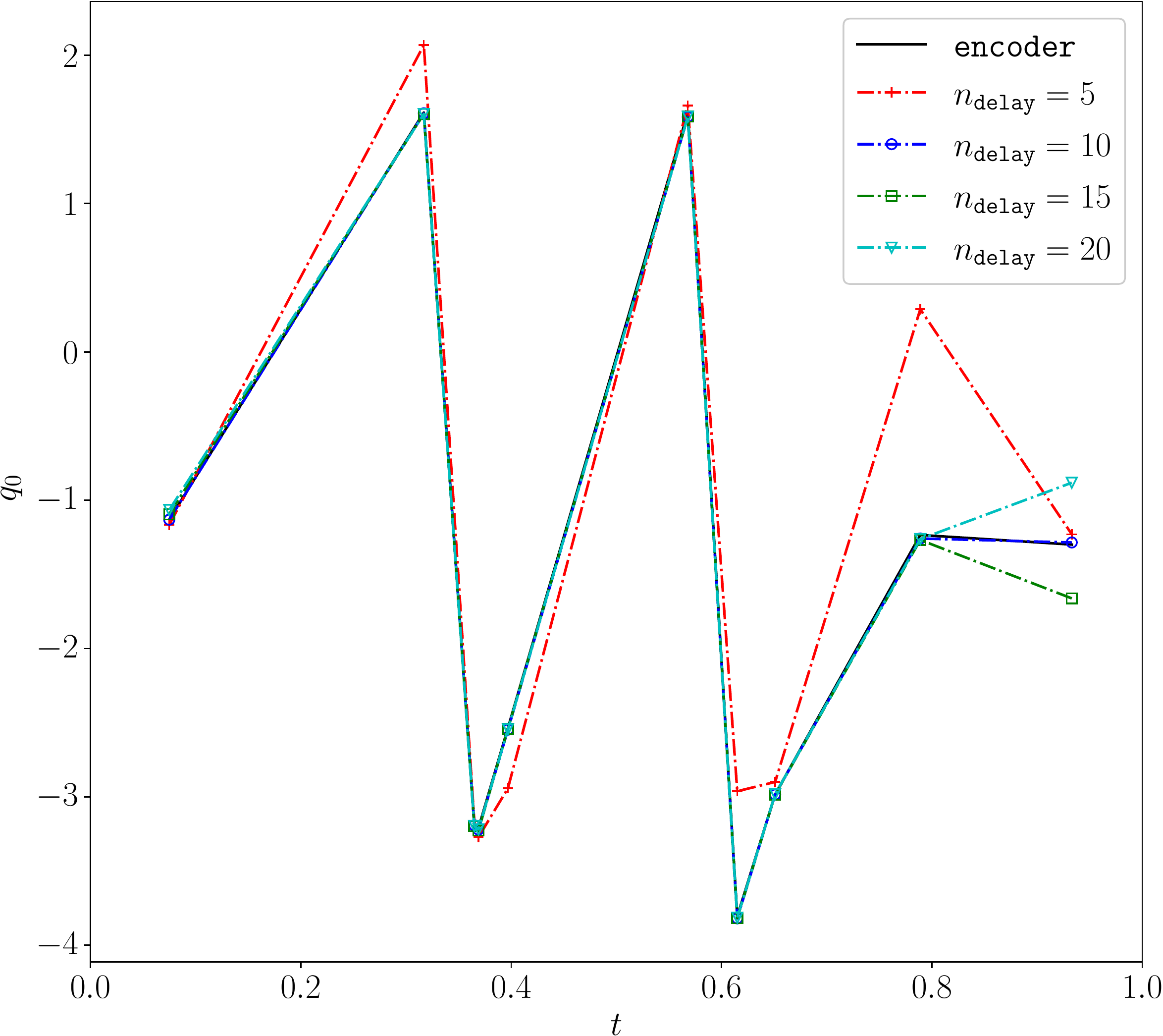}
	\caption{2D Kelvin-Helmholtz instability.
	The temporal evolution of the first latent variable with different $n_{\tt delay}$ for $n_{\tt latent}=8$.}
	\label{fig:2DKHI_q}
\end{figure}

To investigate the dynamics in the latent space,
Figure \ref{fig:2DKHI_err_q_delay} plots the errors $E_{\tt latent}$ with respect to $n_{\tt delay}$ for different $n_{\tt latent}$.
The left figure corresponds to the errors of all the testing times, while the right excludes one trailing time.
Similar to Figure \ref{fig:2DKHI_err_u_CAE_PHODMD}, the errors excluding one trailing time remain almost constant when $n_{\tt delay}$ is large.
Figure \ref{fig:2DKHI_q} shows the evolution of the first latent variable with different $n_{\tt delay}$ for $n_{\tt latent}=8$.
One observes that we can reconstruct accurate dynamics in the latent space.
However, at the ending time, the parametric HODMD tends to be less accurate as $n_{\tt delay}$ increases,
which is due to the formulation of the HODMD.

\section{Conclusion}\label{section:Conc}
This paper has proposed and studied a new non-intrusive reduced-order modeling approach for time-dependent parametrized problems.
Our method is purely data-driven and allows offline and online stages.
Most of the computational costs lie in the offline stage which allows for a fast and accurate prediction for some new times and parameter values at the online stage.
During the offline stage, the deep convolutional autoencoder is first trained.
The encoder reduces the high-dimensional full-order solutions to low-dimensional latent variables,
and the decoder part can recover the full-order solution from the latent space.
The dynamics in the latent space are then modeled by HODMD for each parameter value.
At the online stage, the latent variables at a new time are obtained by the HODMD,
interpolation is utilized to compute the latent variables at a new parameter value,
and the full-order solution is recovered from the latent variables by the decoder.
Numerical tests including the 1D Burgers' equation, 2D Rayleigh-B\'enard convection, and 2D Kelvin-Helmholtz instability
have been conducted to show that our approach can predict the full-order solution at new times and parameter values accurately,
and also works for transport-dominated problems.
In future work, we will consider improving the HODMD method for the latent dynamics.
We will also explore multi-dimensional parameter space and more efficient sampling for the training set.





\end{document}